\documentclass[11pt]{amsart}
\usepackage{amsmath,amsfonts}

\newcommand{\sh}[1]{{#1}}
\newcommand{\rmop}[1]{\mathrm{#1}}

\newcommand{\N}{\mathbb{N}}
\newcommand{\Z}{\mathbb{Z}}
\newcommand{\Q}{\mathbb{Q}}
\newcommand{\R}{\mathbb{R}}
\newcommand{\C}{\mathbb{C}}
\renewcommand{\P}{\mathbb{P}}
\renewcommand{\tt}{\mathfrak{t}}

\renewcommand{\ddot}[1]{{\stackrel{\circ}{#1}}}
\newcommand{\hhat}[1]{\widehat{#1}}

\newcommand{\I}{\mathfrak{I}}
\newcommand{\J}{\mathfrak{J}}
\newcommand{\D}{\mathcal{D}}
\renewcommand{\O}{\mathcal{O}}
\newcommand{\PP}{\mathfrak{P}}

\newcommand{\aarrow}{{\,\to\hspace{-.7em}\to\,}}
\newcommand{\leftmapsto}{{\,\longleftarrow\hspace{-.52em}|\,}}

\newcommand{\Om}{{\Omega}}
\newcommand{\om}{{\varpi}}
\newcommand{\lam}{{\lambda}}
\newcommand{\ep}{{\epsilon}}

\newcommand{\SL}{\mathrm{SL}}
\newcommand{\SU}{\mathrm{SU}}
\newcommand{\U}{\mathrm{U}}
\newcommand{\diag}{\mathrm{diag}}
\newcommand{\Sym}{\mathrm{Sym}}
\newcommand{\Stab}{\mathrm{Stab}}

\newcommand{\ad}{{\mathrm{ad}}}
\newcommand{\Lie}{{\mathrm{Lie}}}
\newcommand{\Ker}{{\mathrm{Ker}}}
\newcommand{\Hom}{{\mathrm{Hom}}}
\newcommand{\Out}{{\mathrm{Out}}}
\newcommand{\Dyn}{{\mathrm{Dyn}}}
\newcommand{\Aut}{{\mathrm{Aut}}}
\newcommand{\id}{{\mathrm{id}}}
\renewcommand{\mod}{\,\mathop{\rm mod}\,}
\newcommand{\Schub}{{\mathfrak{S}}}
\newcommand{\Shat}{\widehat{\Schub}}
\newcommand{\Stil}{\widetilde{\Schub}}
\newcommand{\mult}{\mathrm{mult}}

\newcommand{\cont}{\mathrm{cont}}
\renewcommand{\Im}{\mathrm{Im}}

\newcommand{\proj}{\mathrm{proj}}
\newcommand{\pt}{\mathrm{pt}}
\newcommand{\Span}{\mathrm{Span}}
\newcommand{\ind}{\mathrm{ind}}
\newcommand{\fib}{\mathrm{fib}}
\newcommand{\codim}{\mathrm{codim}}
\newcommand{\eet}{{\perp}}

\newcommand{\What}{\widehat{W}}
\newcommand{\Wtil}{\widetilde{W}}
\newcommand{\Xtil}{\widetilde{X}}

\newcommand{\Gr}{\mathrm{Gr}}
\newcommand{\Grhat}{\widehat{\Gr}}
\newcommand{\Bhat}{\widehat{B}}
\newcommand{\Phat}{\widehat{P}}
\newcommand{\Ghat}{\widehat{G}}

\renewcommand{\t}{(\!(t)\!)}

\newcommand{\xtil}{{\tilde x}}

\begin{document}

\title{Schubert classes of a loop group}
\date{September 2005.  Revised March 2007}
\author[Magyar]{Peter Magyar}

\maketitle
\noindent
In these notes, we
survey the homology of the loop group $\Om K$ of a compact group $K$, also known as the affine Grassmannian $\Grhat$ of a complex loop group
$G[\C\t]$\,.   Using the Bott picture \cite{B} of $H_*(\Om K)$, the homology algebra or Pontryagin ring, we obtain two new results:
\\[-.5em]
\begin{enumerate} \item[A.] 
Factorization of affine Schubert homology classes.
\\[-.5em]
\item[B.] Definition of affine Schubert polynomials representing the affine Schubert homology classes in all types, in terms similar to ordinary Schubert polynomials (Demazure operators on the Borel picture of the cohomology of the flag variety $H^*(K/\,T)$).
\end{enumerate}
\bigskip

The philosophy is that working in the topological category (with compact rather than algebraic groups) gives extra flexibility which simplifies the constructions.  For example, homology factorization happens because the Bott-Samelson variety is topologically a direct product, even though algebraically it is a non-trivial fiber bundle.   Bott and Samelson exploit this flexibility to construct the Pontryagin ring $H_*(\Om K)$ as the symmetric algebra of the homology $H_*(K/\,T)$.  Once this construction is made, the computation of the Schubert classes becomes 
a standard exercise in the theory of Schubert polynomials in
the equivariant cohomology $H^*_T(K/\,T)$.
Thus, one may think of singular homology as a flabby version of equivariant K-theory and representation theory (Demazure modules). 
We work with rational (co)homology throughout this note, but with some care over denominators the results could be refined to integer coefficients.

In this note, 
for concreteness we work with $G=\SL_n\C$,
but unless otherwise noted everything generalizes to a simple, simply connected linear algebraic group of arbitrary type.  In \S1, we survey the relevant background.
We give self-contained statements of Theorem A in \S2.5 and Theorem B 
for $G=\SL_n$ in \S3.6.
In an appendix \S4, we discuss weak order and factorization in the affine Weyl group.

Our formulas are ``global," analogous
to the Schubert polynomials in $H^*(G/B)$ and $H^*_T(G/B)$.  
Let us also mention the unpublished work of Dale Peterson
[P], which contains a different and very powerful point of view on these topics based on the local picture of equivariant cohomology at $T$-fixed points (the nil-Hecke ring of Kostant-Kumar).
Peterson shows that the structure coefficients for $H_*(\Om K)$ are certain Gromov-Witten invariants for the flag variety $G/B$, structure coefficients of the quantum cohomology ring $qH^*(G/B)$.
An alternative proof of Peterson's theorem might
well be possible by combining Bott's realization
of $H_*(\Om K)$ with the known structure of $qH^*(G/B)$.

Using Peterson's picture, Thomas Lam \cite{L} has shown that one version of Lapointe \& Morse's restricted Schur polynomials represent the Schubert classes of $H_*(\Om K)$, so that these polynomials are identical to our Schubert polynomials for $G=\SL_n$\,.  It should be possible to extend this theorem using the computation of affine Schubert polynomials in Theorem B. 

We wish to thank Mark Shimozono, who suggested the questions we examine here.

\section{Background on topology and combinatorics}

\subsection{Loop group and affine Grassmannian}
Let $G\t:=\SL_n\C\t$ be the group of determinant-one matrices with entries in the field of formal Laurent series $\C\t$. 
Then $G\t$ naturally includes the analytic functions 
$f:\C^\times\to G$, 
by identifying a function with its Laurent series.
Let $G[[t]]:=\SL_n\C[[t]]$ be the subgroup of determinant-one matrices with entries in the formal Taylor series $\C[[t]]$.  We may regard $G[[t]]$ as a maximal parabolic subgroup of $G\t$. 
The affine Grassmannian is the quotient space
$
\Grhat:=G\t/\,G[[t]]\,,
$
which we can also realize as the $\C[[t]]$-lattices 
with a fixed virtual dimension in the $\C\t$-vector space $\C\t^n$.

Let $K$ be a maximal compact subgroup of $G$\,, namely  $K=\SU_n=\{A\in M_n(\C)\mid A\bar A^t=I\}$, the group of matrices whose column vectors form an oriented orthonomal basis for $\C^n$ under the Hermitian inner product $\langle (a_1,\ldots,a_n),(b_1,\ldots,b_n)\rangle:=
\sum_i a_i \bar b_i$.  Recall that $K\subset G$ is a maximal compact subgroup, a real Lie group which is not complex algebraic. 

Let $LK$ denote the group of real analytic loops\footnote{
More precisely, we mean the maps which extend to meromorphic functions $f:D\to K$ on the
closed unit disk $D\subset \C$, with poles only allowed at $t=0$.} in $K$, the maps
$f:S^1\to K$, with multiplication performed pointwise: $(f\cdot g)(t):= f(t)\cdot g(t)\in K$. 
Considering $t\in S^1\subset \C^\times$ as a complex parameter, any analytic loop $f(t)\in LK$ has a Laurent series, and we can consider $LK\subset G\t$.
 
Let $\Om K:=\{f(t)\in LK\mid f(1)=1_K\}$ be the normal subgroup of basepoint preserving loops.
We can realize the affine Grassmannian as:
$
\Grhat
\cong \Om K
$\,.
This holds
because Gram-Schmidt orthogonalization allows us to
write any $f(t)\in G\t$ uniquely as $f(t) = f_K(t)\cdot f_P(t)$, where 
$f_K(t)\in \Om K$ and $f_P(t)\in G[[t]]$, so 
$G\t=\Om K\cdot G[[t]]$, a product of transverse subgroups.
%\footnote{
%In fancier language, $G$ is topologically a real vector bundle over its core subgroup $K$, so a loop in $G$ can be retracted to a (based) loop in $K$.
%Analogously in $G$, the compact subgroup $K$ intersects the %Borel subgroup $B$ in the maximal torus $T$, a small subgroup.
%Hence the flag variety $G/B\cong K\!/T$.
%}
Thus the semi-infinite group $\Om K$ provides a canonical slice (set of representatives) for the doubly infinite quotient $\Grhat=G\t/G[[t]]$.
This is what makes the topological theory (using compact groups) simpler than the algebraic theory (using algebraic groups).

It will sometimes be convenient to identify a basepoint-preserving loop in $\Om K$ with the class of its translates by $K$. 
This leads to the realization of the affine Grassmannian as:
$$
\Grhat
\cong 
LK\,/\,(LK\cap G[[t]])=LK/K 
\cong \Om K\,.
$$
where $K\subset LK$ is the group of constant loops $f(t)=k\in K$.
Indeed, we can write any loop in $LK$ as:  
$f(t)=f(t)f(1)^{-1}\cdot f(1)\in\Om K\cdot K$,
and this gives a natural homeomorphism $LK/K\cong\Om K$.
(This is not a group morphism, however.)

\subsection{Pontryagin ring} 

Let us consider the singular homology of $\Grhat=\Om K$.
Although $\Om K$ is infinite-dimensional, 
it is tamely infinite, meaning we can find 
an increasing union of finite-dimensional 
compact Hausdorff subspaces $\Om_1\subset\Om_2\subset
\cdots\subset\Om K$ with $\Om K=\bigcup_{j}\Om_j$.  (We can choose $\Om_j$ from among the Schubert subvarieties we will define below.) \ Thus  $H_*(\Om K)=\displaystyle\lim_\to H_*(\Om_j)$, a direct limit.\footnote{Indeed, we can thicken each $\Om_j$ into a normal neighborhood $U_j\subset\Om K$ which is retractable to $\Om_j$.
Since a singular cycle has compact support, it must live in a sufficiently large $U_j$, so that it is homologous to a cycle in $\Om_j$.}

We have seen that $\Grhat=\Om K$ is topologically a loop group with a continuous (not complex algebraic) multiplication $\mult:\Om K\times \Om K\to\Om K$.
This induces on the singular homology $H_*(\Om K):=H_*(\Om K,\Q)$ a multiplication $\mult_*:H_*(\Om K)\otimes H_*(\Om K)\to H_*(\Om K)$, which together with the usual addition defines the so-called Pontryagin ring.  

It is a basic fact that the Pontryagin multiplication is commutative.  Furthermore, let us consider a different multiplication $\mathrm{conc}:\Om K\times\Om K\to \Om K$ defined by concatenating loops:
$$(f\,\tilde\cdot\,g)(e^{2\pi i x}) := 
\left\{\begin{array}{cl}
f(e^{4\pi i x}) & ,\ \ 0\leq x\leq \frac12\\
g(e^{4\pi i x}) & ,\ \ \frac12\leq x\leq 1\ .
\end{array}\right.
$$
Then the new induced multiplication $\mathrm{conc}_*$ on $H_*(\Om K)$ is the {\it same} as $\mult_*$\,.  
\vspace{-.7em}

{\it Proof:} We construct a homotopy of loops $f\,\tilde\cdot\,g\approx f\cdot g\approx g\,\tilde\cdot\,f$.  For each $-1\leq \epsilon\leq 1$, let $p_\epsilon:[0,1]\to [0,1]\times[0,1]$ be a path in the unit square starting at $(0,0)$ and ending at $(1,1)$, such that $p_{-1}$ goes along the left and top boundaries, $p_0$ goes along the diagonal, and $p_1$ goes along the bottom and right boundaries.  Define $H:[0,1]\times[0,1]\to K$ by $H(x_1,x_2):=f(e^{2\pi ix_1})\,g(e^{2\pi ix_2})$.  Then defining $h_\epsilon:S^1\to K$ by $h_\epsilon(e^{2\pi ix}):=H(p_\epsilon(x))$ gives a continuous family of loops with $h_{-1}=f\,\tilde\cdot\,g$, \ $h_0=f\cdot g$, and
$h_1=g\,\tilde\cdot\,f$. \hfill $\square$
\vspace{-.7em}

As a corollary, the fundamental group of $K$ (or of any Lie group) is abelian, since multiplying in $\pi_1(K)$ means concatenating loops.
Of course, in our case $G=\SL_n\C$ and $K=\SU_n$ are simply connected, so the fundamental group is trivial (and the loop group $\Om K$ is connected).  On the other hand, if we take the quotient of $K$ by the center 
$$
Z=Z(G)=Z(K)=\left\{e^{\frac{2\pi i}{n}j}\, I\ \mid\
 j=0,1,\ldots,n\sh-1\right\}\,,
$$
we get the adjoint groups $G^\ad=G/Z$ and $K^\ad=K/Z$.  Since $G\to G^\ad$ and $K\to K^\ad$ are universal covering maps with fiber $Z$, we have $\pi_1(G^\ad)\cong\pi_1(K^\ad)\cong Z\cong\Z/n\Z$, a cyclic group of order $n$.

\subsection{Coweight lattice and adjoint group}  The subgroup of diagonal matrices $T\subset K$ is a maximal torus: $T\cong (S_1)^{n-1}$.
Letting $\R^n=\R e_1\oplus\cdots\oplus\R e_n$, we may realize 
the Lie algebra $\tt=\Lie(T)$ as:
$$
\tt\cong\{(x_1,\ldots,x_n)\in\R^n\mid
x_1+\cdots+x_n=0\}\,.
$$
The exponential map is a group morphism and universal covering map: $\exp:\tt\to T$, \ $(x_1,\ldots,x_n)\mapsto
\diag(e^{2\pi i x_1},\ldots,e^{2\pi ix_n})$.  
Its kernel is the coroot lattice: 
\begin{eqnarray*}
Q^\vee &:=& \Ker(\exp:\tt\to T)\\
&=&\{\lam\in\Z^n\, 
\mid\,
\lam_1\sh+\cdots\sh+\lam_n=0\}
\ =\ \bigoplus_{i=0}^{n-1}\Z\alpha_i^\vee\,,
\end{eqnarray*}
where $\alpha_i^\vee:=e_i-e_{i+1}$ are the simple coroots.
If we lift the center $Z=Z(G)=Z(K)\subset T$ up to $\tt$, we
get the coweight lattice: 
\begin{eqnarray*}
P^\vee &:=& \exp^{-1}(Z)
\ =\ 
%=\left\{\lam=(\lam_1,\ldots,\lam_n)\in \textstyle\Q^n\left|
%\begin{array}{c}
%\lam_1\sh+\cdots\sh+\lam_n = 0\\
%\lam_i-\lam_{i+1}\in\Z
%\end{array}\right.\right\}
\{ \overline\lam\,\mid\,
\lam\in\Z^n\}
\ =\ \bigoplus_{i=1}^{n-1}\Z\om_i^\vee\,,
\end{eqnarray*}
where we use the projection $\stackrel{-}{\mbox{\ }}\,: \R^n\to\tt$\,, \,
$\overline\lam:=
\lam-\frac{\lam_1+\cdots+\lam_n}{n}(e_1\sh+\cdots\sh+e_n)$\,,
and $\om_i:=\overline{e_1\sh+\cdots+e_i}\in\frac1n\Z^n$ are the fundamental coweights.
Exponentiating $P^\vee$ back into $T$, we see that $Z\cong P^\vee\!/\,Q^\vee$.

Let us identify each vector $\lam\in Q^\vee$ with the straight-line path
$[0,1]\to\tt$, \ $x\mapsto x\lam$ and then exponentiate both sides.  This makes $\lam\in Q^\vee$ correspond to the loop $S^1\to T$, \
$t\mapsto t^\lam:=\diag(t^{\lam_1},\ldots,t^{\lam_n})$, and
in fact $Q^\vee\cong\mathrm{Hom}(S^1,T)$, the homomorphisms from the circle group to the torus.   Similarly, $P^\vee\cong\mathrm{Hom}(S^1,T^\ad)$, where $T^\ad:=T/Z$ is the maximal torus of $K^\ad=K/Z$. 

Consider $\Om K^\ad$, the loops in the adjoint group. 
The connected components of $\Om K^\ad$ are by definition the elements of the fundamental group $\pi_1(K^\ad)$.
Also, the identity component $\Om K^\ad_\circ$ (consisting of homotopically trivial loops) is normal, and the
quotient is the component group $\Om K^\ad\,/\,\Om K^\ad_\circ$ $\cong \pi_1(K^\ad)$.  

We can understand $\Om K^\ad$ purely in terms of $K$.
Every loop $h:[0,1]\to K^\ad$ with $h(0)=h(1)=1$
lifts uniquely to a path $\tilde h:[0,1]\to K$ starting at $\tilde h(0)=1$ and ending at some $\tilde h(1)=z\in Z$.  
The map
$\Om K^\ad\to Z$~, \ $h(x)\mapsto \tilde h(1)=z$ realizes the quotient map to the component group $Z\cong \pi_1(K^\ad)$\,.

A homotopically trivial loop $h(x)\in \Om K^\ad_\circ$
lifts to a loop in $K$, meaning $\tilde h(1)=1$ and
$\tilde h(x)=f(e^{2\pi ix})$ for $f(t)\in\Om K$,
and this identifies $\Om K^\ad_\circ\cong \Om K$.  
As for a component corresponding to
$z=\exp(\mu)\in Z$ with $\mu\in P^\vee$, 
the map $\Om K\to\Om K^\ad$, \ $\tilde h(x)\mapsto \tilde h(x)\cdot \exp(2\pi i x\mu)$ \ (pointwise product) gives a homeomorphism from $\Om K$ to the $z$-component of $\Om K^\ad$.

\subsection{Affine Weyl group}
Let $W=S_n$ be the finite Weyl group,
whose elements $w\in W$ act 
on $\tt\subset\R^n$ by permuting the $n$ basis elements $e_1,\ldots,e_n$, and on $\R^n\aarrow\tt^*$ by permuting the 
dual basis elements $x_1,\ldots,x_n$\,.  
Then $W$ is a Coxeter group whose reflections correspond to the roots of $G$:  for each
pair of a coroot $\alpha^\vee=e_i-e_j\in\tt$ and its corresponding root $\alpha =x_i-x_j\in\tt^*$, we have the reflection
$r_\alpha:\tt\to\tt$, \ $\mu\mapsto \mu-\langle\alpha,\mu\rangle\,\alpha^\vee$.
The simple reflections are $s_1,\ldots,s_{n-1}$ for
$\alpha_i=x_i-x_{i+1}$.

We also have the affine Weyl group $\What=W\times Q^\vee$, a semi-direct product in which $W$ acts on $Q^\vee$.  For $w\in W$, $\lam\in Q^\vee$, an element $w t^{\lam}\in\What$ acts on $\tt\subset\R^n$ by the affine transformation $\mu\mapsto w(\mu+\lam)$, and $wt^\lam=t^{w(\lam)}w$.  Then $\What$ is a Coxeter group whose (affine) reflections
correspond to affine roots, pairs $\hhat\beta=(\alpha,k)$ with $\alpha$ a root of $W$ and $k\in\Z$:
namely, 
$r_{\hhat\beta}=r_{(\alpha,k)}:=r_\alpha t^{k\alpha^\vee}= t^{-k\alpha^\vee}r_\alpha$.

The affine simple roots are
$\hat\alpha_1=(\alpha_1,0)$,
\ldots, $\hat\alpha_{n-1}=(\alpha_{n-1},0)$,
and $\hat\alpha_0:=(-\theta,1)$, where
$\theta=x_1\sh-x_n$ is the highest root and $\theta^\vee=e_1\sh-e_n$ its coroot.\footnote{Warning: For $G$ not simply laced, $\theta^\vee$ is not the highest root of the dual root system.} 
The affine simple reflections are 
$s_i:=r_{\hhat\alpha_i}$ for $i=0,\ldots,n\sh-1$.  Note that
$s_1,\ldots,s_{n-1}\in W\subset\What$,
and $s_0=
r_\theta\, t^{-\theta^\vee}
=t^{\theta^\vee}\! r_\theta$.

We may consider the Weyl group as the quotient $W\cong N_K(T)\,/\,T$, where the normalizer $N_K(T)$ consists of scaled permutation matrices.  We will denote a representative of $w\in W$ by $\ddot w\in N_K(T)$.  Similarly, we have:
$\What\cong N_{LK}(L^{\!0}T)\,/L^{\!0}T$, where
$L^{\!0}T$ is the group of loops $f:S^1\to T$ which extend to analytic functions on the closed unit disk $D\subset\C$.  
Furthermore, the loop $t^\lam\in\Hom(S^1,T)\subset N_{LK}(L^{\!0}T)$ gives a canonical representative for the corresponding element $t^\lam\in\What$, which justifies using the same symbol for both.

There is a technical unpleasantness in this definition of $\What$, however.  The adjoint action of $N_{LK}(L^{\!0}T)$ on $\tt$ has kernel $LT$, not just $L^{\!0}T$, so the adjoint action does not produce the affine action of $\What$ on $\tt$.  This can be remedied by passing to the central extension of $LK$ (or of $G\t$), which is the true affine Kac-Moody group. 

\subsection{Extended affine Weyl group.}
We also have the extended affine Weyl group $\Wtil=W\times P^\vee$ consisting of affine transformations $wt^\lam$ with $w\in W$ and $\lam\in P^\vee$, so that $\Wtil/\What\cong
P^\vee\!/\,Q^\vee\cong Z$.  There is a subgroup $\Sigma\subset\Wtil$ such that $\Wtil= \Sigma\cdot\What$\,,
a semi-direct product in which the subgroup $\Sigma\cong Z$
acts on $\What$.  In fact, $\Sigma=\Stab_{\Wtil}(A_0)$, 
where $A_0$ is the \mbox{fundamental} alcove, a simplex which is a fundamental domain of $\What$ acting on $\tt$: 
$$
A_0:=\{\lam\in\tt\,\mid\, \langle\alpha_i,\lam\rangle\geq 0,\ 
\langle\theta,\lam\rangle\leq 1\}\,.
$$
We can explicitly describe the isomorphism
$\Sigma\cong P^\vee\!/\,Q^\vee$ as follows:
$$\begin{array}{rcl}
\Sigma &\longleftrightarrow& P^\vee/\,Q^\vee\\
\sigma&\longmapsto& \eta=\sigma(0)\\
\sigma=t^{\eta}w_0w_{\eta}&\longleftarrow\hspace{-.5em}\dashv&\lam\,,
\end{array}$$
where $\eta$ is the minimal dominant coweight with $\eta=\lam\mod Q^\vee$; and
$w_0,w_\eta\in W$ are respectively the longest element and the longest element with $w(\eta)=\eta$.
The $\eta$ which occur give canonical representatives for $P^\vee\!/\,Q^\vee$\,:  they are precisely the 
minuscule\footnote{
A coweight $\lam\in P^\vee$ is minuscule if $\langle\alpha,\lam\rangle=0$ or 1 for all positive roots $\alpha$.
}
fundamental coweights $\eta=\om_i^\vee$
(and $\eta=0$ for $\sigma=\id$).    
We can lift $\sigma$ to $LK$ as:
$\ddot\sigma:=t^{\eta} \ddot w_0\ddot w_{\eta}\in LK$.

Every $v\in\Wtil$ has a unique decompostion $v=\sigma_v\hhat v$ with $\sigma_v\in\Sigma$ and $\hhat v\in\What$\,.
To extend the length function to $\Wtil$ we let
$\ell(v):=\ell(\hhat v)$, so that
$\Sigma=\{v\in\Wtil\mid \ell(v)=0\}$.

\subsection{Geometry of the center}

We record some irrelevant facts about the geometry of the center $Z$ and the adjoint group $K^\ad=K/Z$.

Consider $K^\ad$ acting on $K$ by conjugation.  This clearly commutes with multiplication by a central element $z\in Z$ given by: $\tau:K\to K$, \ $k\mapsto z \cdot k$. (This is a deck-shuffling tranformation of the
universal covering map $K\to K^\ad$.) \
Both of these actions are isometries of $K$ (under a suitable Riemannian metric).   

Since every conjugacy class of $K$ intersects $T$ in an orbit of the Weyl group $W$, we can construct the quotient as:
$$
K/K^\ad\cong T/W\cong \tt/(W\times Q^\vee)\cong A_0\,,
$$
where $A_0$ is the fundamental alcove of the previous section.
Now, $\tau=\tau_z$ induces an isometry of the quotient $K/K^\ad$, namely: 
$$\begin{array}{rcl}\overline\tau:A_0&\to& A_0\\
p&\mapsto& p+\eta \ \mod \ (W\times Q^\vee)\,,
\end{array}$$ 
where $\eta=\eta_z$ is the minimal dominant coweight
with $\exp(2\pi i\eta)=z$\,:  this is always a minuscule
fundamental coweight, or $\eta=0$ if $z=1$.
The formula $\sigma=t^\eta w_0w_\eta$ 
lifts $\overline\tau$ to a map $\sigma:\tt\to\tt$ with $\sigma(A_0)=A_0$.

The geometry of $A_0$ is encoded in the affine Dynkin diagram $\widehat\Dyn$:  for example, the walls of $A_0$ correspond to the vertices of $\widehat\Dyn$, the affine wall corresponding to the distinguished 0-node.
A non-trivial isometry $\sigma=\sigma_z$ of $A_0$ defines a graph automorphism of $\widehat\Dyn$ which moves the 0-node to the node $i$ corresponding to the minuscule fundamental coweight $\om_i^\vee=\eta=\sigma(0)$.  Thus we may embed $Z\hookrightarrow \Sym(\widehat\Dyn)$.
Finally, $\widehat\Dyn$ is the Dynkin graph of the Kac-Moody group $G\t$, and any element of $\widehat\Dyn$ induces an outer automorphism of $G\t$, $LK$, and $\Om K$.

To summarize the many guises of the center:
\begin{eqnarray*}
Z&=&Z(G)=Z(K)\\
&\cong& P^\vee\!/\,Q^\vee
=\{0\}\cup\{\eta=\om_i^\vee\ \text{minuscule}\}\\
&\cong& \Wtil/\,\What=\Sigma=\{\sigma\in W\sh\times P^\vee\,\mid\, \sigma(A_0)=A_0\}\\ 
&\cong&\pi_1(K^\ad)
\cong\Om K^\ad/\,\Om K\\
&\hookrightarrow&\text{Isometry}(A_0)\cong\Aut(\widehat{\mathrm{Dyn}})
\cong \Out\,G\t\,.
\end{eqnarray*}

\subsection{Flag varieties}  Let $B\subset G$, the Borel subgroup of upper-triangular matrices in $G$.  
For a coweight $\lam\in P^\vee$, we have the subgroup:
$$
G^\lam:=\Stab_G(t^\lam):=\{g\in G\mid gt^\lam g^{-1}=t^\lam \text{ for all } t\in \C^\times\}\,,
$$
as well as the parabolic subgroup generated by these two subgroups:  $P_\lam:=G_\lam B\subset G$.
The compact version is $K^\lam:=K\cap P^\lam$. 
Note that $T\subset K^\lam$ and $K^{w(\lam)}=wK^\lam w^{-1}$.

We will be concerned with $\om:=\om_1^\vee=\overline{e_1}\in\tt$, the first fundamental coweight, for which we have:
$$P^\om=\{A=(a_{ij})\in \SL_n\C\mid a_{i1}=0 \text{ for } j=2,\ldots,n\}$$ and $K^\om= \SU_n\cap\, (\U_1\times \U_{n-1})$, where 
$\U_1\times \U_{n-1}$ denotes block-diagonal $n\times n$
unitary matrices.  The corresponding partial flag variety is the complex projective $(n\sh-1)$-space: 
$$
X:=G/P^\om=K/K^\om
\cong \P^{n-1}=\C^n\sh-\{0\}\,/\,\C^\times
=S^{2n-1}/S^1\,.
$$
where we consider $S^{2n-1}\subset\R^{2n}\cong\C^n$\,.\,\footnote{
For $n=2$, the circle bundle $S^1\to S^3\to\P^1 \cong S^2$ is known as the Hopf fibration.} Explicitly, we have the homeomorphism:
\begin{eqnarray*}
X=G/P^\om&\stackrel\sim\to& \P^{n-1}=\C^n\sh-\{0\}\,/\,\C^\times\\
A = (a_{ij}) &\mapsto& A\cdot e_1 = [a_{11}:\,\cdots\,:a_{n1}]\,,
\end{eqnarray*}
so that a matrix goes to the projective point
given by its first column.

The full flag variety is the quotient space 
$\Xtil=G/B$, which can be realized as the space
of all nested sequences of complex subspaces
$(V_1\subset V_2\subset\cdots\subset V_n=\C^n)$
with $\dim V_i=i$.
The compact description of this is: $\Xtil\cong K/(K\cap B)=K\!/T$, which is the space of all orthonormal bases of $\C^n$ modulo scaling each basis vector by $S^1$ :  i.e., the flag variety is the space of all splittings of $\C^n$ into orthogonal complex lines.  

There is a right $W$-action on $K\!/T$:
for $x=\ddot x T$ and $w = \ddot w T$, we let $x\cdot w := \ddot x\,\ddot w\,T$.  This action is continuous but not holomorphic.  
For example, if $G=\SL_2\C$ and $w=(12)\in W$, 
we have $\ddot w=${\tiny $
\left[\begin{array}{@{\,}c@{\ }r@{\,}}0&-1\\1&0\end{array}\right]$},
and for a typical flag 
$x=$ {\tiny $\left[\begin{array}{@{\,}c@{\ \,}c@{\,}}1&0\\ a&1\end{array}\right]$}$\in G/B$ we compute:
$$\begin{array}{ccccccc}
G/B &\cong&K\!/T&\stackrel{w}\to&K\!/T&\cong& G/B\\[.3em]
x\mod B&\mapsto& k\mod T&\mapsto& kw\mod T&\mapsto& y\mod B\\[.3em]
\left[\begin{array}{@{\,}rr@{\,}}1&0\\ a&1\end{array}\right]
&\mapsto&
\left[\begin{array}{@{\,}rr@{\,}}1/d&-\bar a/d\\ a/d&1/d\end{array}\right]
&\mapsto&
\left[\begin{array}{@{\,}rr@{\,}}-\bar a/d&-1/d\\ 1/d&-a/d\end{array}\right]
&\mapsto&
\left[\begin{array}{@{\,}rr@{\,}}1\ \ &0\\ 
-1/\bar a&1\end{array}\right]\,,
\end{array}
$$
where $d:=\sqrt{1+a\bar a}$.
Here the first map is Gram-Schmidt orthogonalization, the second map is the $W$-action, and the third map is the holomorphic coordinate function on a chart of $K\!/T\cong G/B\cong\P^1$.
Thus, on this holomorphic chart the map is
$w:a\mapsto -1/\bar a$, which is anti-holomorphic (and orientation-reversing).

\subsection{Classifying spaces}

Consider the direct limit vector space 
$\C^\infty:=\lim_{j\to\infty} \C^j$ with the
$\ell^2$ metric and
its unit sphere $S^\infty$, which
is contractable since its homology is trivial, and which has a free right action of the circle group $S^1$.
The basic example of a classifying space is the quotient $B(S^1)=\P^\infty:=S^\infty\!/S^1$, 
the (tamely) infinite projective space. 
 
Consider $E:=(S^\infty)^n_\ind$, the space of linearly 
independent $n$-tuples of unit vectors
$[\vec v_1,\ldots, \vec v_n]\in(\C^\infty)^n$ thought of as $(\infty\sh\times n)$-matrices\,.
This is contractable with a free right action of the torus $(S^1)^n$ and of the unitary group $U_n$, so the quotients 
$B(S^1)^n=E\!/(S^1)^n$ and $BU_n:=E\!/U_n$ are classifying spaces of the respective groups.   We have $B(S^1)^n\cong(\P^\infty)^n_\ind$, the space of $n$-tuples of linearly independent lines in $\C^\infty$, and
$BU_n\cong\Gr(n,\C^{\infty})$, the Grassmannian of $n$-planes in $\C^\infty$.

For a torus $T\subset (S^1)^n$, we define its classifying space as the free quotient $BT:=E/T$.  For any compact group with an $n$-dimensional representation, $K\subset U_n$, we let $BK:=E/K$.  (The space $E$ serves as a contractable principal bundle for any such $T$ and $K$, so that $E=ET=EK$.)

In our case $T\subset K=SU_n$\,, the space $BT$ will be an $S^1$-bundle over $B(S^1)^n$.  Indeed, if we consider
the dual tautological complex line bundle 
$\D:=\O(1,\ldots,1)$ over $(\P^\infty)^n_\ind$\,, or simply over $(\P^\infty)^n$\,,
we can construct our space $BT$ as the associated $S^1$-bundle 
(Thom bundle) of the complex line bundle $\D$.  
An element of $BT$ is essentially an $n$-tuple
of complex lines $(\C\vec v_1,\ldots,\C\vec v_n)\in
(\P^\infty)^n$, together with
a non-zero volume form\footnote{
More precisely, a $\C$-linear $n$-form on $V$
which agrees in absolute
value with the standard form.}
on the complex $n$-space 
$V=\C\vec v_1\sh+\cdots\sh+\,\C\vec v_n$.
We can also construct $BK$ as the $S^1$-bundle
of the complex line bundle $\D=\O(1)=\wedge^n \mathrm{Taut}^*\to \Gr(n,\C^{\infty})\cong BU_n$\,:
that is, an element of $BK$ is a complex $n$-space $V\subset\C^\infty$ endowed with a non-zero volume form.

We have the $T$-equivariant cohomology ring
$H^*_T(\pt):=H^*(BT)\cong\Sym(\tt^*)$, 
the ring of polynomial functions on $\tt\subset\R^n$.
Indeed, for an element of the coroot lattice 
$\lam\in Q^\vee\subset\tt^*$, we have the
one-dimensional character $\exp(\lam): T\to S^1$
and the associated complex line bundle 
$\mathcal{L}_\lam:=ET\times^T \C_{\exp(\lam)}\to BT$\,.
Then $\lam$ is identified 
with the first Chern class $c_1(\mathcal{L}_\lam)\in H^*(BT)$,
the Poincare dual of the vanishing locus of a section of
$\mathcal{L}$.
In our case we have:
$$
H^*_T(\pt)=
S:=\Q[y_1',\ldots,y_n']/(y_1'+\cdots+y_n')\,,
$$
where $y_1',\ldots,y_n'$ in $(\R^n)^*$ forms a dual basis of
$e_1,\ldots,e_n\in\R^n$.  For historical and combinatorial reasons, we re-index the $y$-variables, letting: $$y_i:=y_{n+1-i}'=y_{w_0(i)}'\,.$$
We define a grading on $\Q[y]$ by writing:
$\dim_\R(y_i)=2$.  
  
The space $BK$ has the same rational cohomology (but \emph{not} 
the same integer cohomology) as the finite quotient $BT\!/W$, where
$W$ permutes the entries of $(\vec v_1,\ldots,\vec v_n)\in
(\P^\infty)^n$:
that is, $H^*(BK)\cong H^*_T(\pt)^W\cong S^W$, the $W$-invariant polynomials in $S$.

\subsection{Cohomology of the flag variety}

For a variety $Y$ endowed with a $T$-action, 
we define the equivariant cohomology ring
$H^*_T(Y):= H^*(Y_T)$, the ordinary cohomology
of the tamely infinite-dimensional induced space: 
$$
Y_T:=ET\stackrel{T}\times Y=
(ET\times Y)\,/\,T\,,
$$
where the free $T$-action on $ET\sh\times Y$ is given by:
$(e,y)\cdot t=(et,t^{-1}y)$\,.

Consider the flag variety $\Xtil=K\!/T$, with
its $T$-equivariant cohomology ring 
$H^*_T(\Xtil):=H^*(\Xtil_T)$.
We have a homeomorphism: 
$$\begin{array}{rcl}
\Xtil_T=ET\stackrel T\times K\!/T
&\stackrel\sim\longrightarrow& 
BT{\displaystyle\mathop{\times}_{BK}}BT\\[1em]
([\vec v_1,\ldots,\vec v_n]\ ,\,k)&\longmapsto&
\left([\vec v_1,\ldots,\vec v_n]\,{\cdot}\, k \ ,\
[\vec v_1,\ldots,\vec v_n]\right)
\ \mod\ T^2\ .
\end{array}$$
From this it is easy to see that
$H^*_T(\Xtil)\cong S\otimes_{S^W} S$.

We can realize the ordinary cohomology $H^*(\Xtil)$
as a quotient of the equivariant cohomology $H^*_T(\Xtil)$ as follows.
Choose an arbitrary basepoint $b_\circ\in BT$, and
consider the fiber bundle:
$$
\Xtil\,\stackrel\fib\to\, ET{\stackrel{T}{\times}} \Xtil\,\stackrel\proj\to\, BT\,,
$$
where $\fib:\Xtil\stackrel\sim\to \Xtil_\circ:=\proj^{-1}(b_\circ)$,
the fiber above the basepoint. 
Then we have a surjective quotient map
$\fib^*:H^*_T(\Xtil)\to H^*(\Xtil)$ with kernel $(\proj^*H^*_T(\pt)_+)$, 
the ideal generated by the pullback classes of positive degree.
Writing this in terms of the polynomial algebra $S$, we get the Borel 
picture of the cohomology:
$$\begin{array}{rcl}
H^*_T(\Xtil)&\aarrow& H^*(\Xtil)\\[.5em]
\displaystyle S\mathop\otimes_{S^W}S
&\aarrow&
{S}/{(S^W_+)}\ .
\end{array}$$

Furthermore, consider the projective space $X=K/K^\om$. The projection $\pi:\Xtil\mapsto X$ induces an inclusion of the equivariant cohomology:
$$\begin{array}{rrcl}
\pi^*:&H^*_T(X)&\hookrightarrow& H^*_T(\Xtil)\\[.5em]
&\displaystyle S^{\om}\mathop\otimes_{S^W}S
&\subset& \displaystyle S\mathop\otimes_{S^W}S\ ,
\end{array}$$
where $S^{\om}$ denotes the invariants under $W^\om$, the Weyl group of $K^{\om}$.  Similarly for the ordinary cohomology groups.  

\subsection{Cohomology generators and relations}
We introduce coordinates for our cohomology rings.
Recall that $BT\to(\P^\infty)^n$ 
possesses $n$ dual tautological bundles $\O(0,\ldots,1,\ldots,0)$.  We let 
$x_i$ (resp.~$y_i'$) denote the Chern class of the  $i^{\mathrm{th}}$ bundle on the first factor 
(resp.~the second factor) of $BT\times_{BK}BT$,
and we recall our convention that 
$y_i:=y_{n+1-i}'\,.$
We now get:
$$
H^*_T(\Xtil)\ \cong\ S{\mathop\otimes_{S^W}} S
\ \cong\ \frac{\Q[x,y]}{\J}\ :=\
\frac{\Q[x_1,\ldots,x_n,y_1,\ldots,y_n]}
{\left(\begin{array}{@{\!}c@{\!}}
x_1+\cdots+x_n=0\\
h(x)=h(y)\ \ \forall\,h\in S^W
\end{array}\right)}\ ,$$
and:
$$H^*(\Xtil)\ \cong\ 
{S}/{(S^W_+)}\ \cong\ 
\displaystyle\frac{\Q[x]}{\J_0}:=
\frac{\Q[x_1,\ldots,x_n]}{(h(x)=0\ \ \forall\,h\in \Q[x]^W_+)}\ .
$$
We also have:
$$H^*_T(X)\ \cong\ 
S^{\om}\otimes_{S^W}S\ \cong\ 
\Q[\xtil,y]/\J^\om\,,$$
where $\xtil:=x_1$ and $\J^\om:=\J\cap S^{\om}$,
and similarly $H^*(X):=\Q[\xtil]/\J_0^\om$\,.
The restriction maps $\fib^*$ are realized by $y=0$\,,
and $\pi^*$ is the obvious inclusion.

We give Grobner bases for our ideals.  (This is one of the few facts whose generalization to arbitrary $K$ is not known.) \ 
Consider the pure lexicographic term-order: 
$$y_1<y_2<\cdots<y_n<x_1<x_2<\cdots<x_n\,.$$
Let $h_d(x_{[i]})$ denote the complete 
symmetric polynomial of degree $d$ in the variables
$x_1,\ldots,x_i$\,, and let $e_d(y)$ denote the elementary
symmetric polynomial of degree $d$ in the variables
$y_1,\ldots,y_n$.
Then the unique reduced Grobner bases of 
our ideals are:
$$\begin{array}{c}
\J=\left(h_1(x_{[n]})\ ,\ e_1(y)\ ,\  
\sum_{i=0}^d\ (-1)^i\ h_{d-i}(x_{[n+1-d]})\
e_i(y)\ ,\ d=1,\ldots,n
\right)\ ,
\\[1em]
\J_0=\left(\, h_{1}(x_{[n]})\ , \ 
h_{2}(x_{[n-1]})\ , \ldots ,\ 
h_{n}(x_{[1]})\, \right)\ ,
\\[1em]
\J^\om=\left(\, e_1(y)\ ,\
\xtil^n-\xtil^{n-1} e_1(y)+\cdots+(-1)^{n-1}\xtil e_{n-1}(y)
+(-1)^n e_n(y)\,\right)\ ,\
\end{array}$$
and $\J^\om_0=(\,\xtil^n\,)$\,, where $\xtil=x_1$.

\subsection{Poincare duality}

It is possible (e.g., by the work of M.E.~Kazarian) 
to define the Poincare duality isomorphism 
$H^{2\ell}_T(\Xtil)\stackrel{\sim}\to H_{\infty-2\ell}^T(\Xtil)$, 
where the right side can be interpreted as the cycles 
of codimension $2\ell$ in the tamely infinite-dimensional 
space $\Xtil_T$.  
Any $T$-invariant subvariety $Y\subset\Xtil$
defines a variety $Y_T\subset\Xtil_T$,
and conversely any finite-codimension subvariety 
$Y'\subset \Xtil_T$ defines a $T$-invariant
subvariety:
$$
Y'_\eet=\eet Y':=\fib^{-1}(Y'\cap\Xtil_\circ)\ \subset\ \Xtil\,. 
$$
We have $Y=(Y_T)_\eet$.  Also
$\codim(Y\,\sh\subset\,\Xtil)=\codim(Y_T\,\sh\subset\,\Xtil_T)$,
as well as $\codim(Y'\,\sh\subset\,\Xtil_T)=\codim(Y'_\eet\,\sh\subset\,\Xtil)$
provided the intersection 
$Y'\cap\Xtil_\circ$ is transversal.

The geometric operation $Y'\mapsto Y'_\eet$ induces
via Poincare duality the cohomology map 
$\fib^*:H^*_T(\Xtil)\to H^*(\Xtil)$\,, 
where $\fib:\Xtil\to\Xtil_T$ 
is the above fiber map:  that is, $\fib^*[Y']=[Y'_\eet]$\,.
For a polynomial cohomology class 
$f(x,y)\in \Q[x,y]/\J\cong H^*_T(\Xtil)$
we have: 
$$
\fib^*f=\eet f(x,y):=f(x,0)\ \in\ 
\Q[x]/\J_0\cong H^*(\Xtil)\,.
$$
The geometric operation $Y\mapsto Y_T$ does not induce
a well-defined map $H^*(\Xtil)\to H^*_T(\Xtil)$, 
since the equivariant class
$[Y_T]$ depends not only upon the homology class $[Y]$,
but upon its $T$-action.

\subsection{Demazure operations}

For a root $\alpha\in\tt^*$\,, 
there are right and left Demazure operations 
$\partial^R_\alpha,\ \partial^L_\alpha$ defined on 
the equivariant cohomology of the flag variety
$H^*_T(\Xtil)$.

We first define the more elementary right operation
$\partial^R_\alpha$ acting on the ordinary cohomology 
$H^*(\Xtil)$.
We have the $\P^1$-bundle $\pi_\alpha:K\!/T\mapsto K\!/K_\alpha$,
inducing on the cohomology $H^*(K\!/T)$ the pullback map $(\pi_\alpha)^*$ and 
the push-forward map\,\footnote{The so-called
Gysin map, realized on DeRham cohomology as
integrating over fibers.}
$(\pi_\alpha)_*$\,.
Then we let:  $\partial^R_\alpha:=(\pi_\alpha)^*(\pi_\alpha)_*:
H^{2\ell}(\Xtil)\to H^{2\ell-2}(\Xtil)\,.$ \

We can interpret $\partial^R_\alpha$ geometrically by
applying Poincare duality to identify $H^{2\ell}(\Xtil)\cong
H_{2L-2\ell}(\Xtil)$, where $2L=\dim_\R(\Xtil)$.
Namely, given a subvariety 
$Y\subset K\!/T=\Xtil$, we lift it to $YT\subset K$,
right multiply by $K_\alpha\supset T$, and push back down to
get the variety $D^R_\alpha(Y):=YK_\alpha/T\subset\Xtil$.
Then we can express the homology operation 
in terms of fundamental classes as:
$$\begin{array}{rcl}
\partial^R_\alpha: H_{2L-2\ell}(\Xtil)&\longrightarrow&
H_{2L-2\ell+2}(\Xtil)\\[.5em]
[Y]&\longmapsto& 
\left\{\begin{array}{cl}
[D^R_\alpha Y]&\text{if}\ \dim_\R(D^R_\alpha Y)=\dim_\R(Y)+2\\
0&\text{otherwise}\,.
\end{array}\right.
\end{array}$$

To define $\partial^R_\alpha$ on $H^*_T(\Xtil)=H^*(\Xtil_T)$,
we can repeat the above word for word, replacing 
$\Xtil=K\!/T$ by $\Xtil_T=K_T\!/T$.  
The equivariant and ordinary 
Demazure operations, both denoted $\partial^R_\alpha$, commute
with the fiber map $\eet :H^*_T(\Xtil)\aarrow H^*(\Xtil)$:
that is, $\eet\circ\partial^R_\alpha=\partial^R_\alpha
\circ\eet$\,.

In terms of the coordinates $\,H^*_T(\Xtil)\,\cong\, \Q[x,y]/\J$\,,
the operator $\partial^R_\alpha$ 
acts as:
$$
\partial^R_\alpha f(x,y):= \frac{f(x,y)- f(r_\alpha x,y)}{x_i-x_j}
\,,$$
where $\alpha=x_i\sh-x_j$ 
and $r_\alpha=(ij)$ switches $x_i$ and $x_j$\,.
Applying $\eet$ to both sides (setting $y=0$) gives the map on $H^*(\Xtil)\cong \Q[x]/\J_0$.

In case $\alpha=\alpha_i$ is a simple root, we write $\partial^R_{(i)}$.
The operators $\partial^R_{(i)}$ satisfy the braid
relations $\partial^R_{(i)}\partial^R_{(i+1)}\partial^R_{(i)}
=\partial^R_{(i+1)}\partial^R_{(i)}\partial^R_{(i+1)}$, etc.,
so we can define $\partial^R_w:=\partial^R_{(i_1)}\cdots
\partial^R_{(i_\ell)}$ for a reduced factorization 
$w=s_{i_1}\!\sh\circ\cdots\!\circ s_{i_\ell}\in W$.
Warning: unless $\alpha$ is a simple root, 
$\partial^R_{\alpha}\neq\partial^R_{r_\alpha}$.

The left operation $\partial^L_\alpha$ acts only on the equivariant cohomology $H^*_T(\Xtil)$.
To define it, first consider the contraction maps
$\kappa_\alpha: K_\alpha{\times^T}\Xtil\to \Xtil$ and
$\kappa_\alpha^T:(K_\alpha{\times^T}\Xtil)_T\to \Xtil_T$,
where $Y_T:=ET \times^T Y$ for any left $T$-space $Y$.
Now let  $\partial^L_\alpha := (\kappa_\alpha^T)_*(\kappa_\alpha^T)^*$\,.
For a geometric intepretation, consider a $T$-invariant subvariety 
$Y\subset\Xtil$, and 
left multiply by $K_\alpha$ to get the
larger $T$-invariant subvariety $D^L_\alpha(Y):=K_\alpha Y\subset\Xtil$.
Then we can express the Demazure operation 
in terms of fundamental classes as:
$$\begin{array}{rcl}
\partial^L_\alpha: H_{\infty-2\ell}^T(\Xtil)&\longrightarrow&
H_{\infty-2\ell+2}^T(\Xtil)\\[.5em]
[Y]_T&\longmapsto& 
\left\{\begin{array}{cl}
[D^L_\alpha Y]_T&\text{if}\ \dim_\R(D^L_\alpha Y)=\dim_\R(Y)+2\\
0&\text{otherwise}\,.
\end{array}\right.
\end{array}$$
On the polynomial ring $\Q[x,y]$, we have:
$$
\partial^L_\alpha f(x,y):= (-1)\,\frac{f(x,y)- f(x,r_\alpha y)}{y_i-y_j}
\,.$$
where $\alpha=y_i\sh-y_j$ and $r_\alpha=(ij)$ switches $y_i$ and $y_j$\,.
(N.B.~the minus in front!) \ 
We once again have $\partial^L_w$, 
and in general $\partial^L_\alpha\neq\partial^L_{r_\alpha}$.

Let us emphasize that in all the above, $\alpha$  
can be any root, not necessarily simple.
We will sometimes write $\partial^x_\alpha$, $\partial^y_\alpha$ instead of $\partial^R_\alpha$, $\partial^L_\alpha$, respectively.

Finally, we consider the Gysin push-forward of the projection $\pi:\Xtil\to X$\,:
$$\begin{array}{r@{\,}ccl}
\pi_*:&H^{2\ell}_T(\Xtil)&\aarrow& H^{2\ell-2M}_T(X)\\[.5em]
&f(x,y)&\longmapsto& \partial^x_{(\om)} f(x,y)\,,
\end{array}$$
where $2M=\dim_\R(\Xtil)-\dim_\R(X)$\,,
and we denote $\partial^x_{(\om)}:=\partial^R_w$, the right
Demazure operator of $w=w^\om$, the longest element of the maximal parabolic $W^\om=\Stab_W(\om)=S_1\times S_{n-1}$\,.
In fact, given $f(x,y)\in H^*_T(\Xtil)$, we may take its Grobner normal form $\bar f(x,y)$ and expand it as:
$$
\bar f(x,y)=\sum_{\!\!\!\!a=(a_2,\ldots,a_n)\!\!\!\!} \bar f_a(\xtil,y)\ x_2^{a_2}\cdots x_n^{a_n}\,,
$$
where $\xtil=x_1$\,.
Then for $a=(n\sh-2,\ldots,2,1,0)$ we have: 
$$
\partial^x_{(\om)}\, f(x,y)=\partial^x_{(\om)}\, \bar f(x,y) =
\bar f_{a}(\xtil,y)\,\mod \J^\om\,. 
$$

\subsection{Schubert polynomials}

The basepoint $\pt := \id\, T\in K\!/T=\Xtil$ corresponds to  $\pt_T\subset\Xtil_T$, which is isomorphic to the 
diagonal $\diag(BT)\subset BT{\times_{BK}}BT$.

The equivariant class of the basepoint,
i.e., the class of the diagonal, corresponds 
via Poincare duality to the double
Schubert polynomial:
$$
\Schub_{w_0}(x,y):=\prod_{i+j\leq n}(x_i-y_j)
=[\pt]_T\ \in\ H^{2L}_T(\Xtil)\,,
$$
where $2L:=\dim_\R(\Xtil)=2\binom{n}{2}$\,.
\\[.5em]
{\it Proof.}  Recall that $BT\times_{BK}BT$ maps to
the space:
$$
B(S^1)^n\mathop{\times}_{BU_n} B(S^1)^n
=\{[v_1,\ldots,v_n]\times[v'_1,\ldots,v'_n]\}\,,
$$
where $v_i,v'_i\in\P^\infty$ and $[v_1,\ldots,v_n]$\,, $[v'_1,\ldots,v'_n]$ run over orthogonal 1-dimensional splittings of some  $n$-dimensional space inside the Hermitean space $\C^\infty$.
The natural coordinates for 
$H^*_T(\Xtil)=H^*(BT\times_{BK}BT)$ are the Chern classes
$x_1,\ldots,x_n$ and $y'_1,\ldots,y'_n$, which can be defined by the loci: $x_i=[v_i\in H]$
and $y'_i=[v'_i\in H]$\,, where $H\subset\C^\infty$ is a hyperplane.

For $i,j\in[1,n]$, consider the locus on which $v_i$ is orthogonal to $v'_j$\,, the divisor: 
$$
D_{ij}:=\left\{\,v_i\cdot v'_j=0\,\right\}\,.
$$
The diagonal $\Delta BT\subset BT\times_{BK}BT$ is the locus:
$$
\Delta BT=\left\{\,v_i=v'_i\,\mid\, i=1,\ldots,n\,\right\}
=\bigcap_{i<j} D_{ij}\,,
$$
a transversal intersection.
Demazure shows that we can express our divisors as:
$
[D_{ij}]=x_i-y'_j\,, 
$
so:
$$
[\pt]_T= [\Delta BT]=\prod_{i<j}(x_i-y'_j)\,.
$$
Substituting $y_i:=y'_{n+1-i}$ gives the desired formula.
 $\square$
\\[-.5em]

For general $w\in W$, the equivariant classes corresponding
to the Schubert varieties
$\Xtil_w:=\mathrm{closure}(B\cdot wT)\subset\Xtil$
are given by the double Schubert polynomials $\Schub_w(x,y)$.  
These can be computed from $\Schub_{w_0}(x,y)$ 
via the right recurrence:
$$
\Schub_{ws_i}(x,y):=\partial^R_{(i)}\Schub_{w}(x,y)
\qquad\text{if}\ \ \ell(ws_i)=\ell(w)\sh-1\,;
$$
or alternatively the left recurrence:
$$
\Schub_{s_iw}(x,y):=\partial^L_{(i)}\Schub_{w}(x,y)
\qquad\text{if}\ \ \ell(s_iw)=\ell(w)\sh-1\,.
$$
Our constructions imply the
non-trivial fact that these two recurrences produce the same
polynomial $\Schub_w$ for each $w$. This fact is equivalent to 
the identity:
$$
\Schub_w(x,y)=(-1)^{\ell(w)}\,\Schub_{w^{-1}}(y,x)\,.
$$

Setting $y=0$ gives us the single Schubert polynomials $\Schub_w(x):=\Schub_w(x,0)$, representing the Schubert
classes in the ordinary cohomology $H^*(\Xtil)$.  These
polynomials can be computed directly from $\Schub_{w_0}(x)=
x_1^{n-1}x_2^{n-2}\cdots x_{n-2}^2 x_{n-1}$ via the right
recurrence:  however, the left recurrence no longer makes sense.
\\[1em]
{\it Example:} For $G=\SL_3\C$\,, we have the equivariant
cohomology ring of the flag variety $\Xtil=\rmop{Flag}(\C^3)$\,: 
$$
H^*_T(\Xtil)\cong\frac{\Q[x,y]}{\J}
\cong \frac{\Q[x_1,x_2,x_3,y_1,y_2,y_3]}
{\left(\begin{array}{c}
x_1\sh+x_2\sh+x_3\ \ ,\ \ y_1\sh+y_2\sh+y_3\\ 
(x_1^2\sh+x_1x_2+x_2^2)
 -(x_1\sh+x_2)\,(y_1\sh+y_2\sh+y_3)
 +(y_1y_2\sh+y_1y_3\sh+y_2y_3)\\
x_1^3-x_1^2\,(y_1\sh+y_2\sh+y_3)+
x_1\,(y_1y_2\sh+y_1y_3\sh+y_2y_3)-y_1y_2y_3
\end{array}\right)}\ ,
$$
and the double Schubert polynomials are:
$$\begin{array}{c}\Schub_{w_0}=(x_1\sh-y_1)(x_1\sh-y_2)(x_2\sh-y_2)\\[.5em]
\Schub_{s_1s_2}=(x_1\sh-y_1)(x_1\sh-y_2)\ ,\ 
\Schub_{s_2s_1}=(x_1\sh-y_1)(x_2\sh-y_1)\\[.5em]
\Schub_{s_1}=x_1\sh-y_1\ ,\ 
\Schub_{s_2}=x_1\sh+x_2\sh-y_1\sh-y_2\ ,\
\Schub_1=1\,.
\end{array}$$
Setting $y_1=y_2=y_3=0$ gives the ordinary cohomology ring:
$$
H^*(\Xtil)\cong\frac{\Q[x]}{\J_0}
\cong\frac{\Q[x_1,x_2,x_3]}
{\left(\begin{array}{@{\,}c@{\,}}
x_1\sh+x_2\sh+x_3\,,\, 
x_1^2\sh+x_1x_2\sh+x_2^2\,,\,
x_1^3
\end{array}\right)}\ ,
$$
with the basis of single Schubert polynomials:
$$\Schub_{w_0}\sh=\,x_1^2x_2\ ,\ 
\Schub_{s_1s_2}\sh=\,x_1^2\ ,\
\Schub_{s_2s_1}\sh=\,x_1x_2\ ,\
\Schub_{s_1}\sh=\,x_1\ ,\
\Schub_{s_2}\sh=\,x_1\sh+x_2\ ,\
\Schub_1\sh=\,1\ .
$$

For the projective space $X=\P^2$, we take $\xtil:=x_1$\,, and we have:
\\[0em]
$$
H^*_T(X)\cong\frac{\Q[\xtil,y]}{\J^\om}
\cong\frac{\Q[\xtil,y_1,y_2,y_3]}
{\left(\begin{array}{@{\,}c@{\,}}
y_1\sh+y_2\sh+y_3\\
\xtil^3-\xtil^2\,(y_1\sh+y_2\sh+y_3)+
\xtil\,(y_1y_2\sh+y_1y_3\sh+y_2y_3)-y_1y_2y_3
\end{array}\right)}
$$
\\[0em]
and $H^*(X)\cong\Q[\xtil]/\J^\om_0
\cong\Q[\xtil]/(\xtil^n)$,
with the Schubert polynomials:
$$\Schub_{s_2s_1}\sh=\,(\xtil\sh-y_1)(\xtil\sh-y_2)\ ,\
\Schub_{s_1}\sh=\,\xtil\sh-y_1\ ,\
\Schub_1\sh=\,1\ .
$$

\subsection{Weyl group actions on cohomology}

The non-holomorphic right action of $W$ on $\Xtil=K/T$\,,
namely $kT\cdot w:=k\ddot wT$\,, induces
a right action on $H^*_T(\Xtil)$.  Recall that each reflection is an orientation-reversing map. 
 
The naive notion of a left $W$-action, $w\cdot kT:=\ddot wkT$\,, is not well-defined.  However, $w\cdot Y:=\ddot wY$ does define an operation on $T$-invariant subvarieties $Y\subset\Xtil$\,,
so we get a left $W$-action on $H^*_T(\Xtil)$.  

Considering $\Xtil_T$ as $BT\times_{BK}BT$\,, we can see these two actions as the natural $W$-actions on the two factors.
That is, we have the natural action of $W$ on $\tt^*$\,,
written $\lam\mapsto w(\lam)$ or in coordinates
$w(x_i)=x_{w(i)}$\,,
and this induces the usual $W$-action on $S=\Sym(\tt^*)$\,.
If a cohomology class is Poincare dual to a subvariety $Y\subset\Xtil$, so that: 
$$
[Y]_T=f(x,y)=f(x,w_0(y'))\in
S\mathop\otimes_{S^W}S\,,
$$
where $y'$ are the natural coordinates and
$y=w_0(y')$\,,  $y_i=y'_{n+1-i}$  according to our convention,
then:
$$
[w\cdot Y]_T=f(x,w_0(wy'))=f(x,w^*(y))
\quad,\quad
[Y\cdot w]_T=(-1)^{\ell(w)}f(w^{-1}(x),y)
\,,
$$
where $w^*:=w_0\,w\,w_0$ is the Weyl duality, and
the sign adjusts orientation.

Now consider the equivariant class of a $T$-fixed point $w=wT\in\Xtil$\,.  If $\pt=\id T$ is the basepoint 
with $[\pt]_T=\Schub_{w_0}=\prod_{i+j\leq n}(x_i-y_j)$\,,
then $w=w\cdot \pt=(-1)^{\ell(w)}\,\pt\cdot w$, so we have:
$$
[w]_T  \ =\ \Schub_{w_0}(x,w^*(y))
\ =\ (-1)^{\ell(w)}\,\Schub_{w_0}(w^{-1}(x),y)
\ \mod\ \J\,.
$$
The two polynomials on the right are not equal,
but equivalent modulo the ideal.

Considering $H^*_T(X)=S^{\om}\otimes_{S^W} S
\cong \Q[\xtil,y]/\J^\om$,
we have $T$-fixed points $u_i = u_iK^\om\in X=K/K^\om$\,,
where $u_i=s_i s_{i-1}\cdots s_2s_1$ for $i=0,\ldots,n\sh-1$
are the minimal coset representatives of $W/W^\om$\,.
Their classes are:
$$
[u_i]_T=\Schub_{u_{n-1}}(\xtil,u_i^*(y))=
\prod_{\!\!\!j\neq n-i\!\!\!}(\xtil-y_j)\,.
$$
Note that $[u_i\sh\in X]_T=\pi_*\,[u_i\sh\in\Xtil]_T$\,.

\section{Affine Schubert varieties and factorization}

\subsection{Minimal parabolic subgroups}  

The elementary unipotent subgroup associated to a  root $\alpha=e_i\sh-e_j$
is $U_{\alpha}:=\{I+xE_{ij}\mid x\in\C\}\cong(\C,+)$, where $E_{ij}$ is an off-diagonal coordinate matrix.  A minimal parabolic is the subgroup generated by a negative elementary subgroup and the Borel:  $P_\alpha=U_{-\alpha} B$, where 
$\alpha$ is a positive root.  
For a simple root $\alpha=\alpha_i=e_i\sh-e_{i+1}$, we write the parabolic as $P_{(i)}$.  We also have the compact form
$K_\alpha:=K\cap P_\alpha$ consisting of matrices which are diagonal except for a copy of $\U_2$ in the $ij$-block, so that $T\subset K_\alpha$.

We make the corresponding definitions for the loop group.
If $\ep_0:G[[t]]\to G$, \ $f(t)\mapsto f(0)$ is the evaluation map, we define the Iwahori sugroup $\Bhat:=\ep_0^{-1}(B)$, which is also a Borel subgroup of $G\t$ thought of as a Kac-Moody group.   For $\alpha=e_i\sh-e_j$ and $k\in\Z$,
we have the affine root subgroup
$U_{(\alpha,k)}:=\{I+xt^kE_{ij}\mid x\in\C\}$.
The positive affine roots are of the form $(\alpha,0)$ and $(\pm\alpha,k)$ for $\alpha$ a positive root of $G$ and $k>0$.
We define the
minimal parabolic $\Phat_{(\alpha,k)}:=U_{(-\alpha,-k)}\Bhat$
whenever $(\alpha,k)$ is positive.
In particular, for the simple roots $\hhat\alpha_i:=(e_i\sh-e_{i+1},0)$ for $i=1,\ldots,n\sh-1$,
and $\hhat\alpha_0:=(e_n\sh-e_1,1)$, we write the minimal parabolic as $\Phat_{(i)}$.  

The compact form of an affine parabolic is:
$K_{(\alpha,k)}:=LK\,\cap\,\Phat_{(\alpha,k)}$. 
Crucially, $K_{(\alpha,k)}$ is once again a finite-dimensional compact Lie group.   
Indeed, $K_{(e_i-e_j,\,k)}$ for $i<j$ is generated by $T\subset LK$ and
a copy of $\SU_2$ in the $ij$-block, in the form:
$$
\left[\begin{array}{@{\,}cc@{\,}}
\!\!\!a & \!\!-\bar b t^k\\
b t^{-k}\!\!& \bar a
\end{array}\right]\,,
$$
where $a,b\in\C$ with $a\bar a+b\bar b=1$.  (This lies in $\SU_2$ for all $t\in S^1$, since $t^{-1}=\bar t$.) \  
Also $K_{(e_j-e_i,\,k)}=K_{(e_i-e_j,-k)}$, the complex conjugate of the above.
 
We have
$K_{(\alpha,0)}=K_{\alpha}\subset K\subset LK$ for all
the finite roots $\alpha$, and
the groups $K_{(0)},K_{(1)},\ldots,K_{(n-1)}$ corresponding to the affine simple roots are all isomorphic to each other via the diagram automorphism of $\widehat\Dyn$.
For example, for $G=\SL_3\C$, the simple-root compact parabolics are:
$$
K_{(0)}=\left[\begin{array}{@{\,}ccc@{\ \ }}
a & 0 & \!\!\!-\bar bt^{-1}\!\!\!\!\!\\
0 & c & \ 0\\
bt&0&\ \bar a
\end{array}\right],\ \
K_{(1)}=\left[\begin{array}{@{\,}ccc@{\,}}
a & \!\!\!-\bar b& 0\\
b&\ \bar a&0\\
0 & \ 0 & c
\end{array}\right],\ \ 
K_{(2)}=\left[\begin{array}{@{\,}ccc@{\,}}
c & 0 & 0\\[.2em]
0 & a & \!\!\!-\bar b\\
0&b&\bar a
\end{array}\right],
$$
where $a,b,c\in\C$ with $(a\bar a+b\bar b)\,c=1$.

\subsection{Schubert varieties in $\Om K$}

We say the factorization 
$v=v_1\cdots v_r\in\What$ is reduced  if 
$\ell(v)=\ell(v_1)+\cdots+\ell(v_r)$, 
and we write this as: 
$v=v_1\!\circ\cdots\circ v_r$.

Let $v=s_{i_1}\!\!\circ\cdots\circ s_{i_\ell}\in\What$ be a reduced factorization into affine simple reflections.  
The spaces we define below will be independent of the choice of factorization for $v$.
Define the algebraic group Schubert variety:
$$
\Ghat_v:=\Phat_{(i_1)}\cdots \Phat_{(i_r)}\subset G\t\ ;
$$
and the compact group Schubert variety:
$$
LK_v:=LK\cap\Ghat_v= K_{(i_1)}\cdots K_{(i_r)}\subset LK\,.
$$
Also define $\Ghat_{\id}:=\Bhat$ and $LK_{\id}:=T$.
For $w\in W$, we will write $K_w\subset K$ instead 
of $LK_w\subset K\subset LK$.

The Schubert variety in the usual sense (inside the affine Grassmannian) is denoted $\Grhat_v$ in the algebraic category and $\Om K_v$ in the topological category:
$$\begin{array}{cccl}
\quad\Grhat_v&:=&(\,\Ghat_v\cdot G[[t]]\,)\,/\,G[[t]]&\subset\ G\t/G[[t]]\ \cong\ \Grhat\\[1em]
\cong\ \Om K_v&:= &
(LK_v\cdot K)\,/\,K &\subset\ LK/K\ \cong\ \Om K\,.
\end{array}$$
Put another way, let $\nu: LK\to\Om K$, \ $f(t)\mapsto f(t)f(1)^{-1}$ be the base-point normalizing map: then $\Om K_v=\nu(LK_v)$.

We may also obtain the Schubert varieties in $\Grhat\cong\Om K$ as the $\Bhat$-orbit closures:
$
\Grhat_{v}=\mathop{\rm closure}\,(\Bhat\cdot \ddot v\subset \Grhat) 
$,
where $\ddot v\in G\t$ is a representative of $v$.  However, this definition has no compact counterpart: the Schubert cell $\Bhat\cdot\ddot v$ is non-compact, so it cannot be the orbit of any compact group.

Note that 
$\Om K_v=\Om K_u$ whenever 
$vW=uW$;  thus the distinct Schubert varieties in
$\Om K$ are indexed by the right cosets
$
\What/W=\{t^\lam W\}_{\lam\in Q^\vee}
$.
Another useful set of coset representatives
is given by the elements of minimal length $m^\lam:=\min(t^\lam W)$\,: we call these the
right-minimal elements of $\What$.
We denote: $$
\Om K_\lam:=\Om K_{t^\lam}=\Om K_{m^\lam}\qquad\text{and}\qquad
LK_\lam:=LK_{m^\lam}\,.
$$
We have a homeomorphism:\footnote{Indeed, the obvious map $LK_\lam/T\to (LK_\lam\cdot K)/K$
is a continuous bijection between compact Hausdorff spaces, and so is a homeomorphism.}
$$
\Om K_\lam\ :=\ (LK_\lam\cdot K)/K\ \cong\ LK_\lam/\,T\,.
$$

Clearly, $LK_v$ and $\Om K_v$ are finite-dimensional, and one can show that $\Grhat_\lam\cong\Om K_\lam$ is a rational complex algebraic variety with
$\dim_\R=2\,\ell(m^\lam)$.  Indeed, one can construct $\Om K$ explicitly as a CW-complex with one even-dimensional cell $\Bhat\cdot t^\lam\cong \C^{\ell(m^\lam)}$ corresponding to each Schubert variety, so that the fundamental cycles form a $\Q$-linear basis of the loop-group homology:
$$
H_*(\Om K)=\bigoplus_{\lam\in Q^\vee} \Q\,[\Om K_\lam]\,.
$$

\subsection{Schubert varieties in $\Om K^\ad$}

We will also need the Schubert subvarieties of $\Om K^\ad$, indexed by $\Wtil/W\cong P^\vee
\cong \Sigma\times Q^\vee$\,, where $\sigma\in\Sigma$ are the elements of length zero.
We say 
$v=v_1\!\circ\cdots\circ v_\ell\in\Wtil$ is a reduced factorization into extended simple roots if we take $v_j\in\{s_0,s_1,\ldots,s_{n-1}\}\cup\Sigma$\,, and the product is reduced (using the length function on $\Wtil$).
For each $\sigma\in\Sigma$ and its representative $\ddot\sigma\in LK$, we let $K_\sigma:=\ddot\sigma T=T\ddot\sigma\subset LK$.
Now for any $v\in\Wtil$, we can define
$LK_v$, \  $\Om K_v$, etc., word-for-word as before, except that now $\Om K_v\subset \Om K^\ad$, etc.
For $v\in\What$, we consider $\Om K_v\subset\Om K\subset\Om K^\ad$.

To any $\lam\in P^\vee$ we associate $\hhat\lam\in Q^\vee$ as follows.  Take
$m^\lam=\sigma_\lam\,\hhat m^\lam$ with $\sigma_\lam\in\Sigma$ and $\hhat m^\lam\in\What$.
Since $\hhat m^\lam$ is again right-minimal, we have 
$\hhat m^\lam=m^{\hhat\lam}$ for some $\hhat\lam$.  Further, let $\eta_\lam=\sigma_\lam(0)$ be the minuscule fundamental weight (or the zero weight) corresponding to $\sigma_\lam$.
Then we have:
$$\Om K_\lam=\Om K_{m^\lam}
=K_{\sigma_\lam}\cdot\Om K_{\hhat m^\lam}
=\sigma_\lam\cdot\Om K_{\hhat\lam}
\subset \Om K^\ad\,.$$

\subsection{Bott-Samelson varieties}

Now suppose $v=v_1 v_2\cdots v_r$ is any factorization
of $v \in \What$ (or $v\in\Wtil$), not necessarily reduced.  
We define the Bott-Samelson variety as the quotient:
$$
\Om K_{v_1,\ldots,v_r}=LK_{v_1}\!{\stackrel{T}\times}
\cdots{\stackrel{T}\times} LK_{v_r}/\,T:=
\frac{LK_{v_1}\!\!\times\cdots\times  LK_{v_r}}
{T\times\cdots\times T}\,
$$
where $LK_{v_1}\!\sh\times \cdots\sh\times LK_{v_{r}}$
has the following right action of $T^r$:
$$
(x_1,\ldots,x_r)\cdot (t_1,\ldots,t_r):=
(x_1t_1\,,\,t_1^{-1}x_2t_2,\cdots,t_{r-1}^{-1}x_rt_r)\,.
$$
This is homeomorphic to the usual algebraic definition of the Bott-Samelson variety (cf.~Demazure).
Indeed $\Om K_{v_1,\ldots,v_r}$ is a complex algebraic variety
with $\dim_\R=2\ell(v_1)+\cdots+2\ell(v_r)$.

Now let $v=v_1\!\circ\cdots\circ v_r$ with
$v$ right-minimal (which implies that $v_r$ is also right-minimal).   Then we have $\Om K_v\cong LK_v/\,T$, and
we can define the contraction map:
$$\begin{array}{rccl}
\cont:&\Om K_{v_1,\ldots,v_r}&\to&\Om K_v\\[.3em]
&(x_1,\ldots,x_r)&\mapsto& x_1\!\cdots x_r\,.
\end{array}$$
This map is continuous and surjective, and it is birational, meaning that $\cont^{-1}(U)\cong U$ for an open neighborhood $U$ of almost any point of $\Om K_v$.
Therefore the induced homology map
$\cont_*:H_*(\Om K_{v_1,\ldots,v_r})\to H_*(\Om K)$ takes the fundamental class
of the Bott-Samelson variety to that of the Schubert variety:
$\cont_*[\Om K_{v_1,\ldots,v_r}]=[\Om K_v]$.

\subsection{Factorization}  

Recall that $m^\lam:=\min(t^\lam W)$\,; that 
$\{m^\lam\}_{\lam\in Q^\vee}$, where $Q^\vee$ is the coroot lattice, are coset representatives for $\What/W$;
and that the Schubert varieties of $\Om K$ are $\Om K_\lam:=\Om K_{m^\lam}=\Om K_{t^\lam}$ for $\lam\in Q^\vee$.

Furthermore, $\{m^\lam\}_{\lam\in P^\vee}$, where $P^\vee$ is the coweight lattice, are coset representatives for $\Wtil/W$, where $\Wtil$ is the extended affine Weyl group.
For $\lam\in P^\vee$, we decompose $m^\lam=\sigma_\lam\hhat m^\lam$ with $\ell(\sigma_\lam)=0$ and $\hhat m^\lam\in\What$, and we define $\hhat\lam\in Q^\vee$ by $m^{\hhat\lam}=\min(\hhat m^\lam W)$.
\\[1em]
{\sc Theorem A}: \ {\it Suppose we have coweights
$\lam,\mu,\nu\in P^\vee$ with 
$\mu\in P^\vee_-$ anti-dominant, such that
$m^\lam\,t^\mu= m^{\nu}\in\Wtil$ and
$\ell(m^\lam)+\ell(t^\mu)=\ell(m^{\nu})$.
Then the Schubert
homology classes in the Pontryagin ring
$H_*(\Om K)$ satisfy:
$$[\Om K_{\hhat\lam}]\cdot [\Om K_{\hhat\mu}]=[\Om K_{\hhat\nu}]\,.$$
}
\\[-.5em]
\noindent
{\it Proof.} For any $\mu\in P^\vee_-$ (meaning $-\mu\in P^\vee_+$ is dominant), 
we have $m^{\mu}=t^{\mu}$ and
$\ell(w t^{\mu}W)\leq \ell(t^{\mu}W)$ for
all $w\in W$.  
Thus
$
K\cdot LK_\mu\cdot K = LK_\mu\cdot K,
$ \
meaning that $LK_\mu\cdot K$ has a left $K$-action. 

Given $m^\lam\circ\, t^{\mu}=m^\nu$ as in the Theorem,
we can write the corresponding Bott-Samelson variety as:
$$
\Om K_{\lam,\mu}:= LK_\lam{\stackrel{T}\times} LK_{\mu}/\,T\cong 
(LK_\lam\cdot K)\stackrel{K}\times (LK_{\mu}\cdot K) / K\,.
$$
Now we can define a homeomorphism by the two inverse maps:
$$\begin{array}{ccc}
\Om K_{\lam,\mu}&\stackrel\sim\longleftrightarrow&\Om K_\lam\times\Om K_{\mu}\\
||&&||\\
(LK_\lam\cdot K){\stackrel{K}\times} (LK_{\mu}\cdot K)/K
&& (LK_\lam\cdot K)/K\times (LK_{\mu}\cdot K)/K\\[1em]
(\,f_1(t)\,,\,f_2(t)\,)&\longmapsto& (\,f_1(t)\,f_1(1)^{-1}\,,\, f_1(1)\,f_2(t)\,)\\[1em]
%(\,f_1(t)\,k\,,\,k^{-1}\,f_2(t)\,k'\,)&\mapsto& (\,f_1(t)\,,\, f_1(0)\,kk^{-1}\,f_2(t)\,k'\,)\\
(\,f_1(t)\,,f_2(t)\,)&\leftmapsto&
(\,f_1(t)\,,\,f_2(t)\,)\ .
%\\[1em]
%(\,f_1(t)\,k\,,\,k^{-1}\,f_1(0)^{-1}\,f_2(t)\,k'\,)&\leftarrow&
%(\,f_1(t)\,k\,,\,f_2(t)\,k'\,)
\end{array}$$
Furthermore, we have the commutative triangle:
$$\begin{array}{c@{\,}c@{\,}l}
\Om K_{\lam,\mu}&\stackrel{\sim}{\longrightarrow}&
\, \Om K_\lam\sh\times\Om K_\mu\\
&\hspace{-1.2em}{}_\cont\!\!\!\!\searrow&\ \ \downarrow
\mbox{\scriptsize$\mult$}\!\!\!\!\\
&&\Om K_\nu\ ,
\end{array}\,,
$$
where $\mult$ is loop group multiplication and $\cont$ is the birational contraction map.
Taking the induced maps on the fundamental
classes in the homology $H_*(\Om K)$
gives:
$$
[\Om K_\lam]\cdot[\Om K_{\mu}]
=\mult_*[\Om K_\lam \sh\times \Om K_{\mu}]
=\cont_*[\Om K_{\lam,\mu}]
=[\Om K_{\nu}]\,.
$$

Now we have 
$[\Om K_\lam]=[\sigma_\lam]\cdot[\Om K_{\hhat\lam}]
$\,, where $[\sigma_\lam]$ is a one-point homology class,
and:
$$
[\sigma_\lam]\cdot[\Om K_{\hhat\lam}]
\cdot[\sigma_\mu]\cdot[\Om K_{\hhat\mu}]=
[\sigma_\nu]\cdot[\Om K_{\hhat\nu}]\,.
$$
But Pontryagin multiplication is commutative and  
$\sigma_\lam \sigma_\mu=\sigma_\nu$\,, so we can cancel
the one-point homology classes from the above equation. 
\hfill $\square$
\vspace{.7em}

We make some comments on the proof.
The key to defining the homeomorphism $\Om K_{\lam,\mu}\stackrel{\sim}{\to}\Om K_\lam\times\Om K_\mu$
is to pick out the $K$-factor of a loop
$f_1(t)\in LK_\lam\cdot K$\ : it is the basepoint $f_1(1)$.
This allows us to disentangle the $K$-twisted product.
Let us see why this proof cannot be adapted to the algebraic category (and so does not prove the analogous results about Demazure modules).  
Algebraically, we can again write the Bott-Samelson variety as:
$$
\Grhat_{\lam,\mu}:=(\Ghat_\lam\cdot G[[t]])\stackrel{G[[t]]}\times
(\Ghat_\mu\cdot G[[t]])\,/\,G[[t]]\,,
$$
which is homeomorphic to $\Om K_{\lam,\mu}$\,.
But now there is no algebraic way to pick out the
 $G[[t]]$-factor of $f_1(t)\in \Ghat_\lam\cdot G[[t]]$\,, so we cannot disentangle the $G[[t]]$-twisted product.  In fact, $\Grhat_{\lam,\mu}$
and $\Grhat_\lam\times\Grhat_\mu$ have non-isomorphic complex structures which are connected by deforming the complex structure on a fixed underlying topological space, as in the recent work of Beilinson-Drinfeld (see \cite{G}).

\section{Affine Schubert polynomials}

\subsection{Homology symmetric algebra}

Let $X=K/K^\om\cong\P^{n-1}$, where $\om:=\om_1^\vee$, the first fundamental coweight.  We have the singular homology: $$H_*(X):=H_*(\P^{n-1},\Q)=\Q 1\oplus\Q h_1\oplus\cdots\oplus\Q h_{n-1}\,,$$
where $1$ is the 0-dimensional homology class of a point and $h_i$ is the $2i$-dimensional homology class of a projective $i$-plane $\P^i\subset \P^{n-1}$.
For our purposes, $H_*(X)$ is a vector space with no natural multiplication.  

Recall that $H_*(X^m)\cong H_*(X)^{\otimes m}$.  Consider
the embedding $\pt:X^m\,\stackrel{\sim}\to\,\pt\times X^m\subset X^{m+1}$\,: by abuse of notation, $\pt$ means either the basepoint 
of $X$ or the inclusion map. 
Define the symmetric power $S^m H_*(X)$ as the quotient of $H_*(X)^{\otimes m}$ by the commutativity relations, and define the symmetric algebra as the direct limit 
$$
S H_*(X)=\lim_{\longrightarrow}
 S^mH_*(X)
$$ 
under the above embeddings, which induce:
$$
\pt^*:S^m H_*(X)\hookrightarrow S^{m+1} H_*(X)\quad,\quad 
y_1\cdot y_2\cdots y_m\mapsto 1\cdot y_1\cdot y_2\cdots y_m\,.
$$
In particular, the unit is $1 = 1\cdot 1 =\cdots$ and any $y\in H_*(X)$ is identified with $y = 1\cdot y= 1\cdot 1\cdot y=\cdots$.  
In our case, we clearly have 
$S H_*(X)\cong \Q[h_1,\cdots,h_{n-1}]\,,$ a polynomial ring.

\subsection{Bott's isomorphism $\Psi$}

Bott defines\footnote{Our mappings $\Psi,\Phi$ are in Bott's notation $g_s,f_s$ respectively.} the mapping:
$$\begin{array}{cccl}
\Psi:&K&\to&\Om K\\
&k&\mapsto& k\,t^\om\,k^{-1}\,t^{-\om}\,.
\end{array}$$
This descends to a map $\Psi:X=K/K^\om\to\Om K$, and it induces a $\Q$-linear map $\Psi_*:H_*(X)\to H_*(\Om K)$ and 
a ring homomorphism $\Psi_{**}:S H_*(X)\to H_*(\Om K)$.
\\[.7em]
{\sc Theorem:} {\it The Bott map $\Psi_{**}$ is an isomorphism from the polynomial ring $\Q[h_1,\ldots,h_{n-1}]=S H_*(\P^{n-1})$ to the Pontryagin ring $H_*(\Om K)$.
This isomorphism preserves the dimension grading, where $\dim(h_i)=2i$.}
\\[.7em]
In the sections that follow, we will paraphrase Bott's proof of this theorem (at least of the surjectivity part).
Our aim is to explicitly describe for each $\lam\in Q^\vee$
the unique polynomial $\Shat_\lam\in SH_*(X)$ with $\Psi_{**}(\Shat_\lam)=[\Om K_\lam]$.
\\[.7em]
{\it Note:} 
For $G$ of other types, the map $\Psi$ is defined in terms of a weight $\om$ with the property that for each simple root $\alpha_1,\alpha_2,\ldots$ of $G$, 
there exists some $w\in W$ such that $\langle\alpha_i,w(\om)\rangle=1$.
(There is always some fundamental coweight 
$\om=\om_i^\vee$ with this property.) \ 
Then we take $X=K/K^\om$ with $H_*(X)\cong
\Q1\oplus\Q h_1\oplus\cdots\oplus\Q h_r$, 
where the basis elements $h_j$ correspond to the 
cosets in the parabolic Bruhat order $W/W^\om$.

The map $\Psi_{**}$ is always surjective, as will be clear
from our arguments below. 
However, $\Psi_{**}$ might not be injective, 
so that $\Shat_\lam=\Psi_{**}^{-1}[\Om K_\lam]$ is no longer a unique polynomial.
Rather, $H^*(\Om K)$ will be isomorphic to 
a quotient of the polynomial ring $SH_*(X)=\Q[h_1,\ldots,h_r]$ by a certain polynomial ideal $\I$, and $\Shat_\lam$ will be well-defined
modulo $\I$.
The algebraic variety $\mathrm{Spec}\, H_*(\Om K)\subset \C^{r}$ defined by the ideal $\I$ is described in Peterson's notes, and is often called the Peterson variety.

\subsection{Bott's map $\Phi$}

We will find it convenient to rephrase Bott's isomorphism
in terms of the $Z$-shifted map: 
$$\begin{array}{r@{\,}c@{\,}c@{\ }l}
\Phi:&K/K^\om&\to&\Om K^\ad\\
&k&\mapsto& k\,t^{\om}\,k^{-1}\,.
\end{array}$$
We have induced maps 
$\Phi_*:H_*(X)\to H_*(\Om K^\ad)$ and
$\Phi_*^m:H_*(X^m)\to H_*(\Om K^\ad)$.
However, since $\Phi_*(1)=t^\om\neq 1\in \Om K^\ad$,
the maps $\Phi_*^m$ do not intertwine with the
inclusions $S^mH_*(X)\hookrightarrow S^{m+1}H_*(X)$, and 
do {\it not} induce a map $\Phi_{**}: 
SH_*(X)\to H_*(\Om K^\ad)$.

Instead, we consider the normalized maps 
$$\hhat{\Phi^m}(x_1,\ldots,x_m):=\Phi(x_1)\cdots\Phi(x_m)
\cdot t^{-m\om}\,,$$
 which do induce a map
$\hhat\Phi_{**}:SH_*(X)\to H_*(\Om K)$.
Because of the commutativity of the Pontryagin product,
we clearly have $\widehat\Phi_{**}=\Psi_{**}$\,.
\\[.7em]
{\sc Lemma:}  {\it Let $\lam\in P^\vee$ correspond to $\hhat\lam\in Q^\vee$, and let $[Y]\in H_*(X^m)$.  Then:
$$
\Phi^m_*[Y]=[\Om K_\lam]\in H_*(\Om K^\ad)
\quad\Longrightarrow\quad
\Psi_{**}[Y]=[\Om K_{\hhat\lam}]\in H_*(\Om K)\,.
$$}
\\[-.8em]
{\it Proof:}  
Let $[t^\lam]$ denote the zero-dimensional 
homology class of the point $t^\lam\in\Om K^\ad$. 
Let $\eta$ be the minimal dominant weight
with $\eta=m\om\mod Q^\vee$: 
since $t^\eta$ and $t^{m\om}$ lie in the same component of 
$\Om K^\ad$, we have
$[t^{\eta}]=[t^{m\om}]$.   
Indeed, $\Im(\Phi^m)$ lies completely in the 
$(m\om)$-component of $\Om K^\ad$, 
and $\Phi^m_*[Y]=[\Om K_\lam]$, so
we have $\lam=m\om\mod Q^\vee$ also, and 
$\Om K_\lam=\Om K_{\hhat\lam}\cdot t^{\eta}$.
  Then:
$$
\Phi^m_*[Y]=[\Om K_\lam]
=[\Om K_{\hhat\lam}]\cdot [t^{\eta}]
=[\Om K_{\hhat\lam}]\cdot[t^{m\om}]\ .
$$
Hence $\Psi_{**}[Y]=\hhat\Phi_{**}[Y]=\widehat{\Phi^m_*}[Y]=[\Om K_{\lam}]\cdot
[t^{-m\om}]=[\Om K_{\hhat\lam}]$, proving the Lemma.

\subsection{Schubert varieties pulled back to $X^m$}  

Given $v=m^\lam\in\What$ 
and a reduced factorization into simple reflections,
we can group the factorization into elements of $W$ separated by the
affine simple reflection $s_0$\,:
$$
v=w_1\circ s_0\circ w_2\circ s_0\circ\cdots\circ w_r\circ s_0\,,
$$
where $w_j\in W$.
Recall that: 
$$
s_0=t^{\theta^\vee\!/2}\,r_\theta\, t^{-\theta^\vee\!/2}
=t^{\om}\,r_\theta\, t^{-\om}\,,
$$
where $\theta=x_1-x_n$ is the highest root,
since $\langle\theta\,,\,\theta^\vee\!/\,2\rangle= 
\langle\theta,\om\rangle=1$\,.  
Expressing this in terms of root subgroups
$K_{(0)}\subset LK$ and $K_\theta\subset K\subset LK$, we have:
$$
K_{(0)}=K_{(-\theta,1)}=t^\om K_{\theta} t^{-\om}\,.
$$
We may write:
$$\begin{array}{rcl}
t^{-\om}=t^{-\pi(e_1)}&=&t^{\pi(e_2+e_3+\cdots+e_n)}\\
&=&t^{u_1(\om)}\,t^{u_2(\om)}\cdots t^{u_{n-1}(\om)}\\
&=&u_1 t^\om u_1^{-1}\ u_2 t^\om u_2^{-1}\,\cdots\,
u_{n-1}t^\om u_{n-1}^{-1}\,,
\end{array}$$
where: $$u_i=s_is_{i-1}\cdots s_1=(i\sh+1,i,\ldots,2,1)\in W=S_n$$
is an $(i\sh+1)$-cycle with $u_i(1)=i\sh+1$.

Now we can express the Schubert class $\Om K_v\subset\Om K$ in terms of
the map $\Phi$ as follows:
$$\begin{array}{rcl}
\Om K_v&=&(K_{w_1} t^\om K_\theta t^{-\om})
\cdots (K_{w_r}\,t^\om\,K_\theta\,t^{-\om})/T\\[.3em]
&=& (K_{w_1}\,t^\om\,K_\theta\,
u_1\,t^{\om}u_1\,u_2 t^{\om}u_2^{-1}\,\cdots\,
u_{n-1}t^{\om}u_{n-1}^{-1})\cdots\\[.3em]
&=& \{k_1\,t^\om\,k_1^{-1}\cdot
(k_1k'_1\,u_1)\,t^{\om}\,(k_1k'_1\,u_1)^{-1}\cdot
(k_1k'_1u_2)\,t^{\om}\,(k_1k'_1u_2)^{-1}\cdots\}\, \\[.3em]
&=&\{\Phi(k_1)\cdot\Phi(k_1k'_1 u_1)\cdots\Phi(k_1k'_1 u_r)\cdot
\Phi(k_1k'_1k_2)\cdot\Phi(k_1k'_1k_2 k'_2 u_1)\cdots\}\,, \\[.3em]
\end{array}$$
where $k_j$ runs over $K_{w_j}$ 
and $k'_j$ runs over $K_\theta$ for $j=1,\ldots,r$.
Hence we define the submanifold $Y\subset X^m=(K/K^\om)^m$,
where $m=nr$, as:
$$
Y := \left\{\left.\begin{array}{@{\!}c@{\ }}
(k_{(1)},\ k'_{(1)}u_1,\,\ldots, \, k'_{(1)}u_{n-1},\\
\ k_{(2)},\ k'_{(2)}u_1,\,\ldots, \, k'_{(2)}u_{n-1},\\
\vdots\\
\ k_{(r)},\ k'_{(r)}u_1,\,\ldots, \, k'_{(r)}u_{n-1})
\end{array}
\right|\begin{array}{@{\,}c@{\!}}
k_{(j)}\!:=k_1k'_1\cdots k_{j-1}k'_{j-1}k_j\ \,,\ \,
k'_{(j)}\!:=k_{(j)}k'_j\\[1em]
k_1\sh\in K_{w_1}\,,\ldots,\, k_r\sh\in K_{w_r}\ \, ,\ \,
k'_1,\ldots,k'_r\sh\in K_\theta
\end{array}
\right\}\,,$$
so that $\Phi^m(Y)=\widehat{\Phi^m}(Y)=\Om K_v=\Om K_\lam$,
and $\Psi_{**}[Y]=[\Om K_\lam]$.
Thus, we have reduced our problem to evaluating the 
fundamental class $[Y]$ as a polynomial in $\Q[h_1,\ldots,h_{n-1}]$. 

We remark that $Y$ is an embedding of a Bott-Samelson 
variety:
$$\begin{array}{rcl}
Y&\cong& K_{w_1}{\stackrel{T}\times} K_\theta
\!\stackrel{T}\times\cdots\stackrel{T}\times\! 
K_{w_r}{\stackrel{T}\times} K_\theta/T 
\ \ \cong\ \ \Om K_{w_1,s_0,\ldots,w_r,s_0}\,.
\end{array}$$
Indeed, we can express the Bott-Samelson variety of
$\Om K$ purely in terms $K_\alpha\subset K$ by virtue of the
isomorphism of $(T\sh\times T)$-spaces $K_{(0)}\cong K_{\theta}$.

\subsection{Demazure operations on multi-projective space}  

We consider the multi-projective space $X^m$, where
$X=K\!/K^\om\cong\P^{n-1}$ and $m=nr$\,.
For a root $\alpha$ of $W$,
we recall the geometric Demazure operation 
$D_\alpha=D^L_\alpha$ which takes
a $T$-invariant submanifold $Z\subset \Xtil^m$
to the new $T$-invariant submanifold:
$$
D_\alpha(Z):=K_\alpha\cdot Z\subset \Xtil^m\,.
$$ 
(Because we will use only left operations, we omit
the superscript $L$.)

For a finite Weyl group element $w=s_{i_1}\!\sh\circ\cdots\!\circ s_{i_\ell}\in W$,
we can define $D_w:=D_{(i_1)}\cdots D_{(i_\ell)}$,
which is independent of the reduced decomposition.
(Once more, $D_\alpha\neq D_{r_\alpha}$.) \ 
Also we recall the base point $\pt=\id\,K^\om\in X$ 
and the $T$-fixed points $u_i=u_i\cdot\pt\in X$ for $i=0,\ldots,n\sh-1$\,. 
Also, we abuse notation so that 
$w^{(k)}\in W^k$ can mean either a $T$-fixed point
of $X^k$ or an inclusion map 
$X^m\stackrel\sim\to w^{(k)}\times X^m\subset X^{k+m}$\,, and similarly with $\pt$\,.

Now let us use these operations to construct the manifold
$Y\subset X^m$ associated to the reduced factorization
$v=w_1\circ s_0\circ\cdots\circ w_r\circ s_0$ in the previous section.  Denote $u^{(n-1)}:=(u_1,\ldots,u_{n-1})\in W^{n-1}$\,. Then we clearly have:
$$
Y=
D_{w_1} \,\pt\, D_\theta
\,u^{(n-1)}\, D_{w_2}\,\pt\, 
D_\theta\,u^{(n-1)}\cdots
D_{w_r} \,\pt\, D_\theta(u^{(n-1)})\,.
$$
We proceed to interpret in polynomial terms each of 
these operations.

Consider the equivariant cohomology $H^*_T(\Xtil^m):= 
H^*(ET{\times^T}\Xtil^m)$.  By the same argument 
as for the case $m=1$, we have:
$$
ET\stackrel T\times\Xtil^m\quad\cong \quad
\overbrace{BT\mathop\times_{BK}\cdots
\mathop\times_{BK} BT}^{m\ \text{factors}}
\mathop\times_{BK}BT
\,.
$$
Taking $H^*_T(\pt)\cong S\cong\Q[y_1,\ldots,y_n]/(y_1\sh+\cdots\sh+y_n\!)$\,,
we get: 
$$\begin{array}{rcl}
H^*_T(\Xtil^m)&\cong&\displaystyle
S\mathop{\otimes}_{S^W} S\mathop{\otimes}_{S^W} \cdots
\mathop{\otimes}_{S^W} S\\[1.5em]
&\cong& 
\displaystyle\frac
{\Q[x^{(1)},x^{(2)},\ldots,x^{(m)},y]}
{
\left(\begin{array}{ @{\,}c@{\!} }
y_1\sh+\cdots+y_n=0\\
h(x^{(1)})=\cdots=h(x^{(m)})=h(y)\ ,\
\forall h\in S^W
\end{array}\right)
}\,,
\end{array}$$
where $x^{(k)}=(x^{(k)}_1,\ldots,x^{(k)}_n)$
and $y=(y_1,\ldots,y_n)$, 
and the dimension grading becomes
$\dim(x^{(k)}_j)=\dim(y_j)=2$.
We also have:
$$\begin{array}{rcl}
H^*_T(X^m)&\cong&\displaystyle
S^{\om}\mathop{\otimes}_{S^W} S^\om\mathop{\otimes}_{S^W} \cdots
\mathop{\otimes}_{S^W} S^\om\mathop{\otimes}_{S^W}S\\[1.5em]
&\cong& 
\displaystyle\frac
{\Q[\xtil_1,\xtil_2,\ldots,\xtil_m,y]}
{
\left(\begin{array}{ @{\,}c@{\!} }
y_1\sh+\cdots+y_n\\[.3em]
\sum_{j=0}^n\, (-1)^{j}\,\xtil_i^{\,j}\,e_{n-j}(y)\ ,\
i=1,\ldots,m\
\end{array}\right)
}\,,
\end{array}$$
where $\xtil_k:=x^{(k)}_1$\,.

The (left) geometric Demazure operation 
$D_\alpha$ induces via Poincare duality 
the cohomology operation:
$$\begin{array}{rcl}
\partial_\alpha^y:H^{2\ell}_T(\Xtil^m)&\to&
H^{2\ell-2}_T(\Xtil_m)\\[.5em]
f(x,y)&\mapsto&\displaystyle
(-1)\,\frac{f(x,y)-f(x,r_\alpha y)}{y_i-y_j}\,,
\end{array}$$ 
where $\alpha=y_i\sh-y_j$\,;
and we again have $\partial_w^y$ for $w\in W$.

Recall the top-degree Schubert polynomial of $H^*_T(X)=S^\om\otimes_{S^W}S$\,, representing the basepoint
of $X$\,:
$$
[\pt\sh\in X]_T=\Schub_{u_{n-1}}(\xtil,y_1,\ldots,y_n)=
\prod_{\ j=1}^{n-1}(\xtil-y_j)\,,
$$
as well as the other $T$-fixed points:
$$
[u_i\in X]_T=\Schub_{u_{n-1}}(\xtil,u_i^*(y))
=\prod_{j\neq i\sh+1} (\xtil-y_j)\,.
$$
The embedding $u^{(n-1)}:X^k\mapsto X^{k+n-1}$ induces the map:
$$\begin{array}{l}
[u^{(n-1)}]_T\, f(\xtil_1,\ldots,\xtil_k,y) := \\[.3em]
\qquad\Schub_{u_{n-1}}(\xtil_1,u_1^*(y))\cdots
\Schub_{u_{n-1}}(\xtil_{n-1},u_{n-1}^*(y))\,
f(\xtil_{n},\ldots,\xtil_{k+n-1},y)\,.
\end{array}$$
We also have the fiber mapping
$\eet:H^*_T(X)\to H^*(X)$\,, \ $f(\xtil,y)\mapsto f(\xtil,0)$\,.  

Now we write:
$$
[Y]=
\eet\,\partial_{w_1}^y \,[\pt]_T\, \partial_{-\theta}^y
\,[u^{(n-1)}]_T\, \partial_{w_2}^y\,[\pt]^T\, 
\partial_{-\theta}\,[u^{(n-1)}]_T\cdots
\partial_{w_r}^y \,[\pt]_T\, \partial_{-\theta}^y\,[u^{(n-1)}]_T\,.
$$
This is a cohomology class: 
$$
[Y]\ \in\ H^{2m(n-1)-2\ell}(X^m)
\ \cong\ 
\left[
\frac{\Q[\xtil_1,\ldots,\xtil_m]}
{\left(\xtil_1^n,\ldots,\xtil_m^n\right)}\right]_{\dim=2m(n-1)-2\ell}
$$ 
where $\ell=\ell(m^\lam):=r+\ell(w_1)+\cdots+\ell(w_r)$
and $\dim(\xtil_j)=2$\,.
This is Poincare dual to the fundmental class:
$$
[Y]\in H_{2\ell}(X^m)\subset \Span_\Q(1,h_1,\ldots,h_{n-1})^{\otimes m}\,.
$$  
Symmetrizing induces the affine Schubert polynomial:
$$
\Schub_\lam=\Sym\,[Y]
\in
\Sym^m \Span_\Q(1,h_1,\ldots,h_{n-1})
=\Q[h_1,\ldots,h_{n-1}]_{\deg\leq m}\,,
$$
This has $\dim(\Schub_\lam)=2\ell$\,, where the dimension grading is defined by $\dim(h_i)=2i$\,; and $\deg(\Schub_\lam)\leq m=nr$\,, where $\deg$ means the usual grading $\deg(h_i)=1$\,.

\subsection{Combinatorics of Schubert classes}

In this section we collect the previous arguments into
a self-contained combinatorial recipe 
to construct our affine Schubert polynomials:
\\[1em]
{\sc Theorem B:} {\it For any $\lam\in Q^\vee$, the polynomial
$$
\Shat_\lam\in \Q[h_1,\ldots,h_{n-1}]= SH_*(X)
$$
defined by the steps (1)--(6) below satisfies
$\Psi_{**}(\Shat_\lam)=[\Om K_\lam]\in H_*(\Om K)$.
} 
\\[-.5em]
\begin{enumerate}

\item We have the affine Weyl group of periodic permutations of $\Z$:
$$\What=\{ \pi\in S_\Z\mid
\pi(i+n)=\pi(i)+n\}\,,$$ which is the semi-direct product of the
finite Weyl group $W=S_n$\,,
and the root lattice: 
$$Q^\vee=\left\{\ t^\lam\ \left|\begin{array}{c}
\lam=(\lam_1,\ldots,\lam_n)\in\Z^n\\
\lam_1+\cdots+\lam_n=0
\end{array}\right.\right\}\,,$$
where $t^\lam(i):=i+n\lam_{(i\mod n)}$\,.
This is a Coxeter group generated by the simple reflections
$s_i=(i,i\sh+1)\in W$ for $i=1,\ldots,n\sh-1$\,,
and by $s_0=r_\theta t^{\theta}$\,, where $r_\theta=(1,n)$
and $\theta=(1,0,\ldots,0,-1)$\,, the highest coroot.

The cosets $\What/W$ have minimal representatives
$m^\lam=\min(t^\lam W)$ for $\lam\in Q^\vee$.
We take a reduced factorization $m^\lam=w_1\sh\circ\,s_0\sh\circ\cdots\sh\circ\,w_r\sh\circ\,s_0$ 
into finite Weyl group elements $w_1,\ldots,w_r\in W$ 
separated by the affine simple reflection $s_0\in\What$.
\\[-.5em]
\item
We let $m=rn$, and we work with $m+n$ variables $\xtil=(\xtil_1,\ldots,\xtil_m)$ 
and $y=(y_1,\ldots,y_n)$\,.
We define the ring:
$$
\frac{\Q[\xtil,y]}{\J^\om}:=
\frac{\Q[\xtil_1,\ldots,\xtil_m,y_1,\ldots,y_n]}
{
\left(\begin{array}{ @{\,}c@{\!} }
y_1\sh+\cdots+y_n=0\\
\sum_{j=0}^n\, (-1)^{n-j}\,\xtil_k^{\,j}\,e_{n-j}(y)\ ,\
k=1,\ldots,m\
\end{array}\right)
}\,,
$$
where $e_i(y)$ denotes an elementary symmetric polynomial.
The given polynomials are a Grobner basis of the ideal
under the term order $y_1<\cdots<y_n<\xtil_1<\cdots<\xtil_m$\,.
\\[-.5em]
\item
For a root $\alpha=y_i-y_j$\,, we have the divided difference operator:
$$
\partial^y_{\alpha}f(\xtil,y):=\frac{f(\xtil,y)- f(\xtil,r_\alpha(y))}{y_i-y_j}\,,
$$
where $r_{\alpha}$ switches the variables $y_i$ and $y_j$\,.
Here $\alpha$ need not be a simple root: for example, it may 
be the longest root $\theta=y_1-y_n$\,.
For $w=s_{i_1}\!\sh\circ\cdots\sh\circ\,s_{i_p}\in W$ a reduced factorization into simple reflections of $W$, 
we define $\partial^y_{w}:= 
\partial^y_{\alpha_{i_1}}\cdots 
\partial^y_{\alpha_{i_p}}$.
\\[-.5em]
\item
For $i=1,\ldots,n$\,, we have the twisted point-class polynomial: 
$$
\PP^{(k)}_i:=
\mathop{\prod_{j=1}}_{j\neq i}^n(\xtil_k-y_j)\,,
$$
and for an interval of $n\sh-1$ integers $$[a\sh+1,a\sh+n\sh-1]=\{a\sh+1,a\sh+2\ldots,a\sh+n\sh-1\}\,,$$
we define:
$$
\PP^{[a+1,a+n-1]}:=
\PP^{(a+1)}_{n-1}\,\PP^{(a+2)}_{n-2}\cdots\,\PP^{(a+n-1)}_{1}\,.
$$
\\[-1.5em]
\item  To the factorization $m^\lam=w_1\sh\circ\, s_0\sh\circ\cdots
\circ\,w_r\sh\circ\, s_0$ we associate a polynomial $\Stil(\xtil,y)\in\Q[\xtil,y]/\J^\om$
by repeatedly applying the divided
difference operators and multiplying by point-class polynomials:
$$
\Stil(\xtil,y):=\pm\,
\partial^y_{w_1}\, \PP^{(1)}_n\,
\partial^y_{\theta}\, \PP^{[2,\,n]}\  \
\partial^y_{w_2}\, \PP^{(n+1)}_n\,
\partial^y_{\theta}\, \PP^{[n+2,\,2n]}\ \cdots
$$
\\[-2.5em]
$$\hspace{5.8em}\cdots\
\partial^y_{w_r}\, \PP^{(rn-n+1)}_n\,
\partial^y_{\theta}\, \PP^{[rn-n+2,\,rn]}\ ,
$$
where the sign 
$\pm=(-1)^{\ell(w_1)+\cdots+\ell(w_r)}$
takes care of orientations.  

We specialize our polynomial to $y=0$\,: 
$$
\Stil(x):=\Stil(x,0)\in
\frac{\Q[\xtil]}{\J^\om_0}
:=\frac{\Q[\xtil_1,\ldots,\xtil_m]}
{(\xtil_1^{\,n},\ldots,\xtil_{m}^{\,n})}\,.
$$
Note that $\deg(\Stil)=m(n\sh-1)-\ell(m^\lam)$\,.
\\[-.5em]
\item We consider $\Q[\xtil]/{\J_0^\om}$ 
as a $\Q$-vector space with a basis of monomials
$\xtil_1^{i_1}\cdots\xtil_m^{i_m}$ with $0\,\sh\leq\, i_k\,\sh\leq\, n\sh-1$\,.
Now we take the $\Q$-linear mapping:
$$\begin{array}{rccl}
\Sym:& \Q[\xtil]/{\tilde\J_0}&\longrightarrow&\Q[h_1,\ldots,h_{n-1}]\\[.5em]
&\xtil_1^{i_1}\cdots\xtil_m^{i_m}&\longmapsto& h_{n-i_1-1}\cdots h_{n-i_m-1}\ ,
\end{array}$$
where $h_0:=1$.
The mapping $\xtil_k^i\mapsto \xtil^i$ 
corresponds to symmetrizing the tensor product, 
and the mapping $\xtil^i\mapsto h_{n\sh-i\sh-1}$ 
undoes Poincare duality, 
turning a cohomology class into a homology class.

Finally, we define:
$$
\Shat_\lam(h):=\Sym\,\Stil(x)\,.
$$
This is a polynomial with degree $\leq rn$ and with
$\dim_\R\Shat_\lam(h)=\ell(m^\lam)$, where
$\dim_\R(h_i):=2i$\,.  It is independent of
the choice of reduced factorization for $m^\lam$.

\end{enumerate}

\subsection{Affine Schubert polynomials for $\SL_3$\,}

Let us first apply the above recipe in the case $n=3$.
We can deduce from Theorem A that all Schubert classes are products of three ``prime" classes:
$$
\Shat_{\rho}\ ,\ \Shat_{-\om_1}\ ,\ \Shat_{-\om_2}\,.
$$
Here we have $\rho:=\om_1\sh+\om_2$\,,
and we write $\om_i$ instead of $\om_i^\vee$ for legibility.
Also, for $\lam\in P^\vee$, we denote 
$\Schub_\lam:=\Schub_{\widehat\lam}$ where
$\widehat\lam\in Q^\vee$ is defined by
$m^{\lam}=\sigma_\lam m^{\widehat\lam}$
with $\ell(\sigma_\lam)=0$\,.
We have: 
$$\Wtil=\langle s_0,s_1,s_2,\sigma\rangle\quad
\text{where}\quad \sigma s_i\sigma^{-1}=s_{(i+1\mod 3)}\,.
$$
To simplify notation further,
we will write $x_k$ instead of $\xtil_k$.
We have:
$$
\PP^{(k)}_i:=\left(x_k\sh-y_1\right)
\left(x_k\sh-y_2\right)\left(x_k\sh-y_3\right)
\ /\ (x_k\sh-y_i)\,.
$$
We compute as follows:
\\[-.2em]
\begin{itemize}
\item $\lambda=\rho=\om_1\sh+\om_2$\,, \ $m^\lambda=s_0$\,, \ $r=\,1$\,, \ $m=\,3$\,.
\\[-.2em]
\begin{itemize}
\item[$\circ$] \ \ 
$
\Stil(x,y):= \PP^{(1)}_3\,\partial^y_{13}\,\PP^{(2)}_2\,\PP^{(3)}_1\\[.3em]
\mbox{\qquad\qquad}=(x_1-y_1)(x_1-y_2)(x_3-y_2)(x_2-y_3)(x_2-y_1)$
\\[-.2em]
\item[$\circ$] \ \ 
$
\Stil(x)=
x_1^2x_2^2x_3$\quad , \quad
$\Shat(h)=h_1$\,.
\end{itemize}
\vspace{1em}

\item $\lambda=-\om_1$\,, \  $m^\lambda=\sigma^{2}s_2s_0$\,, \  $r=1$\,, \  $m=3$\,.
\\[-.2em]
\begin{itemize}
\item[$\circ$] \ \ 
$
\Stil(x,y):= \partial_{32}\PP^{(1)}_3\,\partial^y_{13}\,\PP^{(2)}_2\,\PP^{(3)}_1\\[.3em]
=
(x_1-y_2)(x_2-y_1)
(x_1x_2-x_1x_3+x_2x_3-x_2y_2-x_2y_3+y_2y_3)
$
\\[-.2em]
\item[$\circ$] \ \ 
$
\Stil(x)=
x_1^2x_2^2-x_1^2x_2x_3
+x_2x_2^2x_3$
\\[-.2em]
\item[] \ \ 
$\Shat(h)\ =\ h_2-h_1^2+h_1^2\ =\ h_2$\,.
\end{itemize}
\vspace{1em}

\item $\lambda=-\om_2$\,, \  $m^\lambda=\sigma s_1s_0$\,, \  $r=1$\,, \  $m=3$\,.
\\[-.2em]
\begin{itemize}
\item[$\circ$] \ \ 
$
\Stil(x,y):= \partial_{21}\PP^{(1)}_3\,\partial^y_{13}\,\PP^{(2)}_2\,\PP^{(3)}_1\\[.3em]
=-(x_2-x_3)(x_2-y_3)(x_1-y_1)(x_1-y_2)
$
\\[-.2em]
\item[$\circ$] \ \ 
$
\Stil(x)=x_1^2x_2x_3-x_1^2x_2^2$
\quad,\quad
$\Shat(h)\ =\ h_1^2-h_2$\,.
\end{itemize}
\vspace{1em}

\item $\lambda=-2\om_2$\,, \  
$m^\lambda=\sigma^2 s_0s_2s_1s_0$\,, \  $r=2$\,, \  $m=6$\,.
\\[-.2em]
\begin{itemize}
\item[$\circ$] \ \ 
$
\Stil(x,y):= 
\PP^{(1)}_3\,\partial^y_{13}\,\PP^{(2)}_2\,\PP^{(3)}_1
\partial_{32}\partial_{21}\PP^{(4)}_3\,\partial^y_{13}\,\PP^{(5)}_2\,\PP^{(6)}_1
$
\\[-.2em]
\item[$\circ$] \ \ 
$
\Stil(x)=-x_1^2x_2^2x_3(x_3-x_4)(x_5-x_6)(x_4-x_5)
$
\\[-.2em]
\item[] \ \ 
$\Shat(h)\ =\ (h_1^2-h_2)^2$\,, consistent with Theorem A.
\end{itemize}
\vspace{1em}

\end{itemize}

\subsection{Affine Schubert polynomials for $\SL_4$\,}

For $n=4$, we list the $\lam\in P^\vee$ corresponding to
the ``prime" classes $\Schub_{\widehat\lam}
\in\Q[h_1,h_2,h_3]$.

\begin{itemize}

\item $\lambda=-\om_1$\,, \
$m^\lambda=\sigma^3 s_2s_3s_0$\,, \ 
$\Schub_{\widehat\lam}=h_3$\,.
\\[-.5em]
\item $\lambda=-\om_2$\,, \ 
$m^\lam=\sigma^2 s_0s_3s_1s_0$\,, \ 
$\Schub_\lam=h_2^2-h_1h_3$\,.
\\[-.5em]
\item $\lambda=-\om_3$\,, \ 
$m^\lam=\sigma s_2s_1s_0$\,, \
$\Schub_\lam=h_1^3-2h_1h_2+h_3$\,.
\\[-.5em]
\item $\lambda=\rho=\om_1\sh+\om_2\sh+\om_3$\,, \ 
$m^\lam=\sigma^2s_2s_1s_3s_0$\,, \ 
$\Schub_\lam=h_2(h_1^2-h_2)$\,.
\\[-.5em]
\item $m^\lam=\sigma^2s_1s_3s_0$\,, \
$\Schub_\lam=h_1h_2-h_3$\,.
\\[-.5em]
\item $m^\lam=\sigma^2 s_3s_0$\,, \ 
$\Schub_\lam=h_2$\,.
\\[-.5em]
\item $m^\lam=\sigma^2 s_1s_0$\,, \ 
$\Schub_\lam=h_1^2-h_2$\,.
\\[-.5em]
\item $m^\lam=\sigma^2 s_0$\,, \ 
$\Schub_\lam=h_1$\,.

\end{itemize}

Notice the bit of mysterious extra factorization for $\Schub_\rho$\,.

\subsection{Calculations for higher $\SL_n$\,}

For general $n$, the ``prime" Schubert classes include
$\Schub_{-\om_i}$, the classes corresponding to the 
fundamental anti-dominant translations $t^{-\om_i}=m^{-\om_i}$.
(We again identify $\om_i$ with $\om_i^\vee$, etc.) \ 
We can decompose this as:
$$
t^{-\om_i}
=(s_is_{i-1}\cdots s_2s_1\sigma)^{n-i}
=w_iw_0\sigma^{-i}\,,
$$
where $w_i=\max W^{\om_i}$, the longest element of a maximal parabolic, and $\sigma$ is the automorphism of the fundamental alcove with:
$$
\sigma s_i\sigma^{-1}=s_{i+1}\quad,\quad\text{where}\  s_n:=s_0\ .
$$

Now, $\lam=-\om_i\not\in Q^\vee$, so $\Schub_\lam$ means
$\Schub_{\widehat\lam}$\,, where $\widehat\lam\in Q^\vee$ is defined by
$m^{\widehat\lam}=\sigma^i t^{-\om_i}$\,.
Indeed, Shimozono has calculated that, 
letting $k:=\min(i,n\sh-i)$, we have: 
$$
\widehat\lam=k\alpha_i+\sum_{j=1}^{k-1}(k\sh-j)(\alpha_{i-j}+\alpha_{i+j})\,.
$$
For example for $n=6$, we have: 
$$\begin{array}{ccc}
\widehat{-\om_1}=\alpha_1&
\widehat{-\om_2}=\alpha_1+2\alpha_2+\alpha_3\\[-.1em]
&&\widehat{-\om_3}=\alpha_1+2\alpha_2+3\alpha_3+2\alpha_4+\alpha_5\ .
\\[-.5em]
\widehat{-\om_5}=\alpha_5&
\widehat{-\om_4}=\alpha_3+2\alpha_4+\alpha_5
\end{array}$$

The other prime classes are $\Schub_\lam$ for $\lam\in Q^\vee$
with $(m^\lam)^{-1}$ lying inside the fundamental box,
as described in the Appendix.
The fundamental box has maximal alcove 
$$
m^\rho=\ t^\rho w_0\in\Wtil\,,
$$
where $\rho=\sum_{i=1}^{n-1}\om_i$\,, which is also the half-sum of the positive coroots. 
This element is an involution: 
$$
(m^\rho)^{-1}=w_0t^\rho=t^{-w_0(\rho)}w_0 =m^\rho\,.
$$
We can decompose this as $m^\rho=\iota\,\hat m^\rho$,
where $\hat m^\rho\in\What$ and:
$$
\iota = \left\{\begin{array}{c@{\,}l}
1 & \text{ for $n$ odd}\\
\sigma^{n/2} &\text{ for $n$ even\,.}
\end{array}\right.
$$
The alcoves in the fundamental box are the right weak order interval $[1,\hat m]_R$\,, and the corresponding prime classes are: 
$$
[1,\hat m^{-1}]_L=[1,\iota\,\hat m\,\iota]_L
=\left\{w\in\What\mid \,\iota\,\hat m\,\iota = u\circ w\text{ for some } u
\right\}\,.
$$

\section{Appendix: Weak order on an affine Weyl group}

\newcommand{\hW}{\widehat{W}}
\newcommand{\tW}{\widetilde{W}}
\newcommand{\hatt}{{ {}^\wedge\hspace{-.1em} t }}
\renewcommand{\v}{^\vee}
\newcommand{\omv}{\om\v}
\newcommand{\al}{\alpha}
\newcommand{\alv}{\al\v}
\newcommand{\sig}{\sigma}
\newcommand{\Desc}{\mathrm{Desc}}
\newcommand{\Inv}{\mathrm{Inv}}

\newcommand{\leqB}{\stackrel{\mathrm B}{\leq}}
\newcommand{\leqL}{\stackrel{\mathrm L}{\leq}}
\newcommand{\leqR}{\stackrel{\mathrm R}{\leq}}

\newcommand{\lhat}[1]{{ {}^\wedge\!#1 }}
\newcommand{\rhat}[1]{{ #1^{\!\wedge} }}

This note addresses the following question concerning optimal factorization in an affine Weyl group.

\subsection{Question}  Given an affine Weyl group element $x\in\hW$, minimal in its coset $xW$ with respect to the finite Weyl group $W$, find the maximal dominant coweight $\lam$ such that there exists a factorization $x=y\cdot t_{-\lam}$ with $\ell(x)=\ell(y)+\ell(t_{-\lam})$, where $y\in \tW$, the extended affine Weyl group, and $t_{-\lam}\in\tW$ is an anti-dominant translation.

\subsection{Answer}  
Let $\omv_1,\ldots,\omv_n$ be the fundamental coweights of $W$.  For a vector $\gamma = c_1\omv_1+\cdots c_n\omv_n$ with $c_i\in\R$, define its floor as the integral coweight:
$$
\lfloor\gamma\rfloor:=\lfloor c_1\rfloor\omv_1+\cdots+\lfloor c_n\rfloor\omv_n\,,
$$
where $\lfloor c\rfloor$ denotes the integer floor of a real number.  

Choose a tiny positive vector $\ep:=\ep_1\omv_1+\cdots+\ep_n\omv_n$, where $\ep_i>0$ are sufficiently small positive real numbers. 
Then the dominant coweight asked for is:
$\lam=\lfloor x^{-1}(\ep)\rfloor$, which is independent of the choice of $\ep$.

\subsection{Examples}
\begin{itemize}

\item[a.]
For $x=w t_\mu$, where $w\in W$ and $\mu$ is a dominant coweight, this gives:
$$
\lam=\mu-\!\!\!\!\sum_{i\in D(w^{-1})}\!\! \omv_i\,,
$$
where $D(w^{-1}):=\{ i\in[1,n]\ \mid\, w^{-1}(\al_i)<0\}$, the descent set of $w^{-1}$.

\item[b.]
Let $\hW$ be of type $A_2^{(1)}$ and $x=s_1s_0$.  We have $s_1(\omv_1)=-\omv_1+\omv_2$, \ $s_2(\omv_2)=\omv_1-\omv_2$, and $s_i(\omv_j)=\omv_j$ for $1\leq i\neq j\leq n$.   Also, $s_0=t_{\theta\v}r_{\theta}$, where $\theta\v=\omv_1+\omv_2$ and $r_{\theta}(\omv_1)=-\omv_2$, \ $r_{\theta}(\omv_2)=-\omv_1$.
Then $\ep=\ep_1\omv_1+\ep_2\omv_2$, and $x^{-1}(\ep)=(1\sh-\ep_1\sh-\ep_2)\omv_1+(1\sh+\ep_2)\omv_2$, so
$\lam=0\,\omv_1+1\,\omv_2=\omv_2$.

\item[c.]
Let $\hW$ be of type $C_2^{(1)}$ with Dynkin diagram
$0\sh\Rightarrow1\sh\Leftarrow\, 2$, and $x=s_0s_1s_0$. 
We have $s_1(\omv_1)=-\omv_1+2\omv_2$, \ $s_2(\omv_2)=\omv_1-\omv_2$, and $s_i(\omv_j)=\omv_j$ for $1\leq i\neq j\leq n$.   Also, $s_0=t_{\theta\v}r_{\theta}$, where $\theta\v=\omv_1$ and $r_{\theta}(\omv_1)=-\omv_1$, \ $r_{\theta}(\omv_2)=-\omv_1+\omv_2$.
Then $\ep=\ep_1\omv_1+\ep_2\omv_2$, and $x^{-1}(\ep)=\ep_1\omv_1+(2\sh-2\ep_1\sh-\ep_2)\omv_2$, so
$\lam=0\,\omv_1+1\,\omv_2=\omv_2$.

\end{itemize}

\subsection{Root systems and alcoves.}  
We proceed to prove the above answer, as well as record some general facts for reference(cf. \cite{H}, \cite{S}).

Let $\R^n$ and $(\R^n)^*$ be dual vector spaces.  
A positive definite inner product allows one to identify $\R^n$ and $(\R^n)^*$, but we will use only the natural pairing $\langle\ ,\ \rangle$ between $(\R^n)^*$ and $\R^n$.
We define a root system by a set
$\Delta\subset(\R^n)^*$ of root vectors and a set
$\Delta\v\subset\R^n$ of coroot vectors.  We have the bases of simple roots
$\al_1,\ldots,\al_n\in(\R^n)^*$ and simple coroots $\alv_1,\ldots\alv_n\in\R^n$, together with their dual bases, the fundamental weights $\om_1,\ldots,\om_n\in(\R^n)^*$ and fundamental coweights $\omv_1,\ldots,\omv_n\in\R^n$.
For $\alpha\in\Delta$, we have the reflection $r_\alpha:\R^n\to \R^n$ , $\mu\mapsto\mu - \langle\alpha,\mu\rangle\,\alv$.  The simple reflections $s_1,\ldots,s_n$ generate the finite Coxeter group $W$.

For a vector $\lam\in\R^n$ we have the translation $t_\lam:\R^n\to\R^n$ , \mbox{$\mu\mapsto\mu+\lam$,} so that $w\,t_\lam\, w^{-1}=t_{w(\lam)}$ for $w\in W$. 
We define the affine Weyl group $\hW = W\propto t_{Q\v}$ and the extended affine Weyl group $\tW = W\propto t_{P\v}$,
where $Q\v = \Z\alv_1\oplus\cdots\oplus\Z\alv_n$ is the coroot lattice and $P\v=\Z\omv_1\oplus\cdots\oplus\Z\omv_n$ is the coweight lattice.
As a Coxeter group, $\hW$ has generators $s_0,s_1,\ldots,s_n$, where $s_0=t_{\theta\v}\,r_{\theta}=r_{\theta} t_{-\theta\v}$. 

The affine reflections $r_{\al,k}\in\hW$ correspond to the hyperplanes $$
H_{\al,k}=\left\{ \mu \mid \, \langle\al,\mu\rangle = k\right\}
$$
for $\al\in\Delta_+$ and $k\in\Z$.
The (closures of) connected components of $\R^n\setminus\{H_{\al,k}\}$ are called {\it alcoves}: they are simplicies provided $\Delta$ is an irreducible root system.  We define the fundamental alcove 
$$\begin{array}{rcl}
A_0&:=&
\{\mu\mid \, 0\leq \langle\al,\mu\rangle\leq 1\ \text{for all}\ \al\in\Delta_+\}\\[.3em]
&=&\{\mu\mid \, \langle\al_1,\mu\rangle,\ldots,\langle\al_n,\mu\rangle\geq 0\ ,\ \langle\theta,\mu\rangle\leq 1\}\,.
\end{array}$$
This is a fundamental domain of $\hW$ acting on $\R^n$, and each alcove $A$ corresponds to a unique $w\in\hW$ with $A = wA_0$.  Henceforth we identify $w\in\hW$ with its alcove $wA_0$.  

Since the neighbors of $A_0$ are $s_iA_0$, any two neighboring alcoves $A,A'$ are of the form $A=wA_0$ and $A'=ws_iA_0$ for some $i=0,\ldots,n$; and we label the common facet of $A$ and $A'$ with $i$.  
The length function on $\hW$ (or on the set of alcoves) is defined as: 
$
\ell(w)= \#\ \text{hyperplanes separating $w$ from $1$}\,.
$

The extended $\tW$ is the symmetry group of the affine hyperplane arrangement $\{H_{\al,k}\}$, or of the set of alcoves, and the quotient $\tW/\hW$ is isomorphic to the symmetry group of $A_0$ itself: $\sig A_0=A_0$ for each minimal coset representative $\sig\in\tW/\hW$.  That is, $\tW=\Sigma\propto\hW$,
where $\Sigma:=\{\sig\in\tW\mid \sig(A_0)=A_0\}$\,. 
Every element $x\in\tW$ can be uniquely factored as: 
$$
x=\lhat x\,\sigma=\sigma\,\rhat x
$$
where $\sigma\in\Sigma$ and
$\lhat x,\rhat x\in\hW$.  
%Indeed, $\lhat x=\sigma\,\rhat x\,\sigma^{-1}$.
We have $xA_0=\lhat xA_0$,
and we define $\ell(x):=\ell(\lhat x)=\ell(\rhat x)$.
In particular,  for each $\lam\in P\v$ we have $t_\lam=\hatt_\lam\,\sigma$ with $\hatt_\lam\in\hW$ and some $\sigma\in \tW/\hW$.

We also have the fundamental chamber 
$C_0=\{\mu\mid\,
\langle \al_1,\mu \rangle,\ldots,
\langle \al_n,\mu \rangle \geq 0 
\}$, 
and a translation $t_\lam$ lies in $C_0$ exactly when $\lam$ is dominant.  Furthermore, for $\lam$ a dominant coweight, $wt_{\lam}$ lies in the chamber $wC_0$, and $t_{-\lam}$ lies in $-C_0 = w_0C_0$.

Lusztig has given the following length formula for 
$wt_\lam\in\tW$, where $w\in W$ and $\lam\in P^\vee$:
$$
\ell(wt_\lam)=\sum_{\alpha\in\Delta_+}|\langle\lambda,\alpha\rangle|+\ell(w)\,.
$$

\subsection{Strong and weak orders.}
There are three important orders on $\hW$ and $\tW$.  The {\it strong} or {\it Bruhat order} $\leqB$ has covering relations $x\leqB r_{\al,k}x$, where $r_{\al,k}$ is any reflection with $\ell(r_{\al,k}x)=\ell(x)+1$.  We have $x\leqB y$ iff $x^{-1}\leqB y^{-1}$. 

In this note, we will mainly consider the weak orders:
the {\it left weak order} $\leqL$ has covers $x\leqL s_ix$, where $s_i$ ($i=0,\ldots,n$) is a simple reflection with $\ell(s_ix)=\ell(x)+1$.  The analogous {\it right weak order} $\leqR$ has covers $x\leqR xs_i$.  
Since $x\leqL y$ iff $x^{-1}\leqR y^{-1}$, the left and right orders are combinatorially equivalent.  
However, the right order is easier to picture in terms of alcoves: the covering relation $x\leqR xs_i$ steps from an alcove $A=xA_0$ to its neighbor $A'=xs_iA_0$ through a facet labelled $i$, provided the facet separates $A'$ from the fundamental $A_0$.

These orders become very simple on dominant translations $t_\lam$.  That is, let $\lam,\mu$ be dominant coweights. 
Then we have: 
$$
t_{\lam}\leqB t_{\mu} \quad\Longleftrightarrow\quad
\mu\sh-\lam \in Q\v_+:=\N\alv_1\oplus\cdots\oplus\N\alv_n
\quad\stackrel{\text{def}}\Longleftrightarrow\quad
\lam\stackrel{Q}\leq \mu\,,
$$
$$
\hatt_{\lam}\leqR \hatt_{\mu} \quad\Longleftrightarrow\quad
\mu\sh-\lam \in P\v_+:=\N\omv_1\oplus\cdots\oplus\N\omv_n
\quad\stackrel{\text{def}}\Longleftrightarrow\quad
\lam\stackrel{P}\leq \mu\,.
$$
In the second line we must use $\hatt$ instead of $t$, since we can only have $t_\mu\leqR t_\lam$ when $t_\mu\equiv t_\lam\!\mod \hW$, meaning $\mu-\lam\in Q\v$.

When restricted to an orbit $W\cdot t^{\lam}$ for $\lam\in P^\vee_+$, the affine Bruhat order becomes
the finite Bruhat order on the parabolic quotient $W\!/\,W_\lam$. If $W_\lam=\{1\}$, we have for $w,y\in W$\,:\\[-.5em]
$$wt_\lam\leqB yt_\lam\quad\Longleftrightarrow\quad w\leqB y
\quad\Longrightarrow\quad
w(\lam)\stackrel{Q}\geq y(\lam)\,.$$
Warning:  The converse of the last implication does \emph{not} hold.\footnote{
Example: $W=S_{n+1}$\,, \ $\lam=\rho=(n\sh+1,n,\cdots,2,1)\in P^\vee_+\subset\R^{n+1}$. We have:
\\[-1em]
$$
w\leqB y\ \ \Longleftrightarrow\ \ 
\{w(1),\ldots,w(i)\}\leq \{y(1),\ldots,y(i)\}\ \ \text{for}\ \ 
i=1,\ldots,n\,,
$$
\\[-2.3em]
$$
w(\rho)\stackrel{Q}\geq y(\rho)\ \ \Longleftrightarrow\ \ 
w(1)\sh+\cdots\sh+w(i)\leq y(1)\sh+\cdots\sh+y(i)\ \ \text{for}\ \ i=1,\ldots,n\,.
$$
For $w=2314$\,, \ $y=4123$\,, the second condition holds but not the first.
}
To characterize the weak order on $W$, define
the inversion set:
$$
\Inv(w):=\Delta_+\cap w^{-1}\Delta_-
=\{\alpha^{(1)},\ldots,\alpha^{(n)}\}\,,
$$
where $\alpha^{(k)}:=
s_{i_n}s_{i_{n-1}}\cdots s_{i_{k+1}}(\alpha_{i_k})$\,.
Now we have: 
$y\leqL w$ iff $\Inv(y)\subset\Inv(w)$\,, \ 
and $y\leqR w$ iff $\Inv(y^{-1})\subset\Inv(w^{-1})$\,. 
To define the strong and left weak orders on $W\!/\,W_\lam$, we replace the cosets $wW_\lam$ by their shortest (or longest) representatives and take the orders induced from $W$.

The orders can be adapted to parabolic quotients of the affine Weyl group.  Consider $\hW/W=\{t_\lam W\}_{\lam\in Q\v}$, the right cosets of the affine Weyl group by the finite Weyl group.   Pictorially, each coset corresponds to a translate $t_\lam P_0$ of the Weyl polytope $P_0:=\bigcup_{w\in W} wA_0$ by an element of the coroot lattice.   Let $\min(xW)$, $\max(xW)$ denote the minimal and maximal element in the coset.
We define: $\ell(xW)=\ell(\min(xW))$, and 
$xW\leqL s_ixW$ whenever $\ell(s_ixW)=\ell(xW)+1$.
In fact, 
$$xW\leqL yW
\quad\Longleftrightarrow\quad
\min(xW)\leqL\min(yW)
\quad\Longleftrightarrow\quad
\max(xW)\leqL\max(yW)\,.
$$

The left cosets $Wx\in W\backslash\hW$
can best be pictured by their minimal representatives, the alcoves inside the fundamental chamber: $A\subset C_0$ with $A=\min(Wx)=t_\lam w$, where $\lam\in Q\v_+$ and $w\in W$
with $w$ minimal in the coset $w\,\mathrm{Stab}_W(\lam)$.
We define $Wx\leqR Wy$ so that: $xW\leqL yW$ iff $Wx^{-1}\leqR Wy^{-1}$.  

The weak orders extend in an obvious way to $\tW/W$ and $W\backslash\tW$.  The Bruhat order on $\tW/W$ can be defined similarly, and I believe the minimal covers of $xW=vt_\lam W$ for $\lam\in P\v_+$ and $v\in W$ are:
$$
vt_\lam W\ \leqB\ \left\{\begin{array}{cl}
 r_\al v\, t_\lam W &\text{if }\ \ell(r_\al v)=\ell(v)+1\\
 r_\al v\, t_{\lam-\alv_i}W &\text{if }\ \ell(r_\al v)=\ell(v)-1\,.\\
\end{array}\right.
$$
I do not know a good criterion for when $vt_\lam\leqB wt_\mu$.

\subsection{Factorization and weak order.}
A {\it reduced factorization} is an expression $x=y\cdot z$, where $x=yz\in\tW$ and $\ell(x)=\ell(y)+\ell(z)$.  This is closely related to the weak order:  for $x,y,z\in\hW$,
$$
x =y\cdot z 
\quad\Longleftrightarrow\quad
x\leqL z
\quad\Longleftrightarrow\quad
x\leqR y\,.
$$ 
For $x,y,z\in\tW$, this becomes:
$$
x =y\cdot z 
\quad\Longleftrightarrow\quad
\rhat x\leqL \rhat z
\quad\Longleftrightarrow\quad
\lhat x\leqR \lhat y\,,
$$ 
since if $x=y\cdot z$ with $y=\lhat y\sig$, \ $z=\lhat z\sig'$, then  
$x=\lhat x\sig\sig'$ and
$\lhat x = \lhat y \cdot (\sig\, \lhat z\, \sig^{-1})$ with
$\sig\,\lhat z\, \sig^{-1}\in\hW$.  
The converse means: if $\lhat x\leqR\lhat y$, then 
for any $\sig,\sig''\in\Sigma$, there exists $z\in\tW$ such that  $\lhat x\sig''=
\lhat y\sig\cdot z$.
Similarly for $\rhat x$.  

Now we can rephrase the Question.  Given $x\in \hW$, minimal in its coset $xW$, we want the maximal $\lam\in P\v_+$ such that for some $y\in\tW$, 
we have: $x=y \cdot t_{-\lam}$, meaning 
$x^{-1}=t_{\lam}\cdot y^{-1}$, 
meaning $\hatt_{\lam}\leqR x^{-1}$.

\subsection{Question reformulated:}  Given $x^{-1}\in\hW$, minimal in its coset $Wx^{-1}$, find the $\stackrel P\leq$-maximal coweight $\lam\in P\v_+$ such that $\hatt_{\lam}\leqR x^{-1}$.
\\[.5em]
This has an obvious pictorial answer.  The group $t_{P\v}$ acts on $\R^n$ with fundamental domain equal to the parallelopiped spanned by the fundamental coweights: 
$$
B_0:=
\{\mu\mid 0\leq\langle\al_i,\mu\rangle\leq 1\text{ for }
i=1,\ldots,n\}\,
$$
We call $B_0$ the {\it fundamental box}, and its translates
$B_\lam:=t_\lam B_0$ cover $\R^n=\bigcup_{\lam\in P\v} B_\lam$.
Each box is a union of $N$ alcoves, where
$N = |W|\, /\, [\tW:\hW] = \mbox{$|W|\,/\, [P\v\!:Q\v]$}$.
Furthermore, the fundamamental chamber is the union of the dominant boxes: \ 
$C_0=\bigcup_{\lam\in P\v_+} B_\lam$.

Now we have: $\hatt_\lam\leqR y$ iff $y A_0\subset B_\lam$.
(The $\Rightarrow$ is easy. The $\Leftarrow$ is because \ldots) \ 
Hence, given $y=x^{-1}\in\hW$ 
we have:
$$
\max\{ \mu \mid \hatt_{\mu}\leqR y\}=\lam\quad\text{where}\quad
yA_0\subset B_\lam\,.
$$
Now, if we take any point 
$
\ep=\ep_1\omv_1+\cdots+\ep_n\omv_n
$ 
in the interior of $A_0$,
we have: 
$$
yA_0\subset B_\lam
\quad\Longleftrightarrow\quad
y(\ep)\in B_\lam
\quad\Longleftrightarrow\quad
\forall\, i :
\langle\al_i,\lam\rangle<
\langle\al_i,y(\ep)\rangle<
\langle\al_i,\lam\rangle+1\,.
$$
Therefore, the answer to our question is exactly:
$$
\lam:=
\lfloor\langle\al_1, y(\ep)\rangle\rfloor\,\omv_1+\cdots
+\lfloor\langle\al_n, y(\ep)\rangle\rfloor\,\omv_n\,,
$$
which was to be shown.

\vspace{2em}
\noindent 
Peter Magyar\\
\verb|magyar@math.msu.edu|\\
Dept of Math, Wells Hall\\
Michigan State University\\
East Lansing, MI 48823\\

\end{document}